\documentclass[twoside,a4paper,reqno,11pt]{amsart}

\usepackage[top=28mm,right=30mm,bottom=28mm,left=30mm]{geometry}
\usepackage{amsfonts, amsmath, amssymb, mathrsfs, bm, latexsym, stmaryrd, array, mathtools, bold-extra, xcolor}
\usepackage[colorlinks,linkcolor=cyan,anchorcolor=blue,citecolor=red,backref=page]{hyperref}

\usepackage{caption}
\usepackage{listings}
\usepackage{arydshln}
\usepackage[all]{xy}
\usepackage{multirow}
\usepackage{enumitem}
\setlist[enumerate]{font={\rm},itemsep=0.2\baselineskip}
\setlist[enumerate,1]{label={(\roman*)}}
\setlist[enumerate,2]{label={(\arabic*)}}


\newtheorem{thm}{Theorem}

\newtheorem{theorem}{Theorem}[section]
\newtheorem{proposition}[theorem]{Proposition}

\newtheorem{corollary}[theorem]{Corollary}

\newtheorem{lemma}[theorem]{Lemma}

\theoremstyle{definition}
\newtheorem{definition}[theorem]{Definition}

\newtheorem{example}[theorem]{Example}

\newtheorem*{conj*}{Conjecture}

\theoremstyle{remark}

\def\GL{{\mathrm{GL}}}
\def\GammaL{{\Gamma \mathrm{L}}}

\def\AGammaL{{\mathrm{A}\Gamma \mathrm{L}}}
\def\SL{{\mathrm{SL}}}
\def\PSL{{\mathrm{PSL}}}

\def\Sp{{\mathrm{Sp}}}
\def\PSp{{\mathrm{PSp}}}
\def\SU{{\mathrm{SU}}}

\def\CGammaSp{{\mathrm{C}\Gamma\mathrm{Sp}}}

\def\A{{\mathrm{A}}}

\def\F{{\mathrm{F}}}
\def\E{{\mathrm{E}}}

\def\calB{{\mathcal{B}}}

\def\calH{{\mathcal{H}}}

\def\bbZ{{\mathbb{Z}}}
\def\bbF{{\mathbb{F}}}

\def\Tr{{\mathrm{Tr}}}

\def\bfZ{{\mathbf{Z}}}

\def\Aut{{\mathrm{Aut}}}
\def\Out{{\mathrm{Out}}}
\def\End{{\mathrm{End}}}
\def\Gal{{\mathrm{Gal}}}

\def\Hol{{\mathrm{Hol}}}

\def\leqs{\leqslant}
\def\geqs{\geqslant}

\def\Magma{{\sc Magma}}
\def\Sz{{\mathrm{Sz}}}

\begin{document}

\title{The finite groups with three automorphism orbits}

\author{Cai Heng Li}

\address{SUSTech International Center for Mathematics, and Department of Mathematics, Southern University of Science and Technology, Shenzhen, Guangdong, China}
\email{lich@sustech.edu.cn {\text{\rm(Li)}}}

\author{Yan Zhou Zhu}

\address{SUSTech International Center for Mathematics, and Department of Mathematics, Southern University of Science and Technology, Shenzhen, Guangdong, China}
\email{zhuyz@mail.sustech.edu.cn {\text{\rm(Zhu)}}}

\begin{abstract}
A complete classification is given of finite groups whose elements are partitioned into three orbits by the automorphism groups, solving the long-standing classification problem initiated by G.\,Higman in 1963.
As a consequence, a classification is obtained for finite permutation groups of rank $3$ which are holomorphs of groups.
\end{abstract}


\maketitle

\section{Introduction}

For a group $N$, the natural action of the automorphism group $\Aut(N)$ on $N$ partitions the elements of $N$ into orbits, called the \textit{automorphism orbits} of $N$.
We remark that an automorphism orbit is also called a \textit{fusion class}.
The study of the relations between groups and their automorphism orbits was initiated by G.\,Higman~\cite{higman1963Suzuki} in 1963 when he studied finite $2$-groups whose involutions form an automorphism orbit, and since then automorphism orbits of groups have received considerable attention from group theorists.

For convenience, a group $N$ is called \textit{$r$-orbit} if $\Aut(N)$ partitions the elements of $N$ into $r$ orbits.
Since elements in the same automorphism orbit have the same order, the trivial group is the only $1$-orbit group, and elementary abelian $p$-groups with $p$ prime are the only finite $2$-orbit groups.
It is worth mentioning that there are examples of non-abelian infinite $2$-orbit groups, see~\cite[Corollary 1.3]{osin2010Small}.
The class of finite $3$-orbit groups contains some interesting families of groups, and the problem of classifying such groups has received considerable attention in the literature.
The main purpose of this paper is to classify all finite $3$-orbit groups.

Obviously, a $3$-orbit group has order divisible by at most two primes, and so it is either a $p$-group or a $\{p,q\}$-group, where $p,q$ are distinct primes.
It is not hard to see that finite $3$-orbit abelian groups are homocyclic groups with exponent $p^2$ with $p$ prime.
Finite $\{p,q\}$-groups that are $3$-orbit groups have been classified in~\cite{devillers2014normal,dornhoff1970imprimitive,laffey1986Automorphism,maurer1997Groups} (also see Example~\ref{exam:pq}), which are proved to be certain Frobenius groups.
M{\"a}urer and Stroppel~\cite{maurer1997Groups} study groups $3$-orbit groups under additional structural assumptions.
Recently, finite $2$-groups that are $3$-orbit groups have been classified independently by Bors and Glasby~\cite{bors2024Finite}, and the authors~\cite{li2024proof}, with two typical examples being Sylow $2$-subgroups of unitary groups $\SU(3,2^f)$ and Suzuki groups $\Sz(2^{2f+1})$ for $f\geqs 1$ (see Section~\ref{sec:exam}).
The problem of classifying $3$-orbit groups is therefore reduced to the the case of $p$-groups with $p$ an odd prime.

As comparison, we would like to mention some interesting research about automorphism orbits.
Studies of $r$-orbit groups for small values of $r>3$ mainly focus on non-solvable groups.
Specifically, $\A_5$ is the only $4$-orbit non-solvable group in~\cite[Theorem 3]{laffey1986Automorphism}, and $\PSL(2,p)$ for $p\in\{7,8,9\}$ are the only $5$-orbit non-abelian simple groups~\cite{stroppel1999finite}.
The finite simple $r$-orbit groups for $r\leqs 100$ are determined in~\cite{jafari2021automorphism}, and the non-solvable $r$-orbit groups for $r=5, 6$ are classified in~\cite{bastos2016finite,dantas2017Finite}.
Finite groups with some special automorphism orbits also receive a lot of attention.
Shult proved that finite $p$-group for odd prime $p$ in which the elements of order $p$ form an automorphism orbit is abelian~\cite{shult1969finite}.
In 1976, Gross studied finite $2$-groups whose involutions form an automorphism orbit~\cite{gross1976automorphic}, and the authors have recently completed the classification of such groups recently~\cite{li2024proof}.
Influenced by Gross' work, finite groups whose elements of the same order lie in the same automorphism orbit, so-called \textit{AT-groups}, were studied by Zhang~\cite{zhang1992finite}.
Praeger and the first-named author~\cite{li1997Finitea} gave a description of finite groups in which elements of the same order are fused or inverse-fused.
Recently, Bors, Giudici and Praeger studied the relation between the number of automorphism orbits and the size of the spectrum of a non-solvable group~\cite{bors2023Automorphisma}.

Recall that a non-abelian $p$-group $N$ is said to be a \textit{special $p$-group} if $\bfZ(N)=\Phi(N)=N'$ is elementary abelian.
In this paper, for a prime power $q=p^f$, we denote by $q^{m+n}$ a non-abelian special $p$-group $N$ with $|\Phi(N)|=q^m$ and $|N/\Phi(N)|=q^n$.
We remark that finite non-abelian $3$-orbit $p$-groups are special $p$-groups since each of them has a unique non-trivial proper characteristic subgroup.
For odd prime $p$, finite non-abelian $3$-orbit $p$-groups are special $p$-groups of exponent $p$.
It is known that the $n$-generator free group in the variety of groups having exponent $p$ and nilpotency class $2$ is finite, which is denoted by $\calH_{n,p}$ in this paper (see Definition~\ref{def:hei}).
Hence, for odd prime $p$, each non-abelian $n$-generator $3$-orbit $p$-group is isomorphic to a quotient group $\calH_{n,p}/W$ for some $W<\bfZ(\calH_{n,p})$.
To classify $3$-orbit groups, we have to determine all ``desired'' subgroups $W$.
Denote by $\Lambda^2_{\F}(V)$ the \textit{exterior square} of vector space $V$ over $\F$.
We remark that $\bfZ(\calH_{n,p})$ can be ``viewed'' as the exterior square of $\calH_{n,p}/\bfZ(\calH_{n,p})$ in some sense (see Lemma~\ref{lem:oddp} for more details).
The following theorem, an analog of \cite[Theorem 1.2]{li2024proof} which is for $p=2$, describes certain representations of finite transitive linear groups of characteristic $p$.

\begin{thm}\label{thm:transmod}
	Let $V=\bbF_p^n$ for odd prime $p$, $G\leqs \GL(V)$ be non-solvable, and let $W<\Lambda^2_{\bbF_p}(V)$.
	If $G$ acts transitively on non-zero vectors of both $V$ and $\Lambda^2_{\bbF_p}(V)/W$, then one of the following statements holds.
	\begin{enumerate}
		\item $\dim \Lambda^2_{\bbF_p}(V)/W=1$.
		\item $G^{(\infty)}\cong\SL(3,p^{n/3})$ and $\dim \Lambda^2_{\bbF_p}(V)/W=\dim V=n$.
		\item $G^{(\infty)}\cong\Sp(2m,p^{n/(2m)})$ and $\Lambda^2_{\bbF_p}(V)/W$ is a trivial $\bbF_pG^{(\infty)}$-module with dimension at most $n/(2m)$.
	\end{enumerate}
\end{thm}

This result plays an important role in proving the second main theorem of the paper, which presents a classification of finite $3$-orbit groups, stated as follows.

\begin{thm}\label{thm:3orbits}
	A finite group is a $3$-orbit group if and only if it is one of the groups listed in Table~$\ref{tab:3aogp}$, where $p,q$ are distinct primes, and $m,n$ are positive integers.
\end{thm}

\begin{table}[ht!]
	\centering
	\caption{Finite $3$-orbit groups}
	\label{tab:3aogp}
	\begin{tabular}{lllll}
    \hline
	&$N$& $\Aut(N)$ & Conditions & Ref \\
	\hline
	1&$\bbZ_{p}^{n(q-1)}{:}\bbZ_q$&$\AGammaL(n,p^{q-1})$&Example~\ref{exam:pq}&\ref{lem:pq}\\
	2&$\bbZ_{p^2}^n$&$\GL(n,\bbZ/p^2\bbZ)\cong p^{n^2}.\GL(n,p)$&&\ref{lem:ab}\\
	3&$A_2(n,\theta)\cong 2^{n+n}$& $2^{n^2}{:}\GammaL(1,2^n)$&$|\theta|\neq 1$ is odd&\ref{thm:2group}\\
	4&$\SU(3,2^n)_2\cong 2^{n+2n}$ & $2^{2n^2}{:}\GammaL(1,2^{2n})$ & &\ref{thm:2group}\\
	5&$P(\epsilon)\cong 2^{3+6}$ & $2^{18}{:}(\bbZ_7{:}\bbZ_9)$ &&\ref{thm:2group}\\
	6&$A_p(n,\theta)\cong p^{n+n}$ & $p^{n^2}{:}\GammaL(3,p^{n/3})$& $p$ is odd and $|\theta|=3$ &\ref{lem:autoapsl}\\
	7&$q^{1+2m}_{+}\cong p^{n+2mn}$ & $p^{2mn^2}{:}(\Sp(2m,q){:}\GammaL(1,q))$& $q=p^n$ and $p$ is odd &\ref{lem:nmpncgamma}\\
	8&$q^{1+2m}_{+}/U\cong p^{n_0+2mn}$ & $p^{2mn_0n}{:}(\Sp(2m,q){:}\GammaL(1,q)_U)$&Example~\ref{exam:quot} &\ref{lem:quotaut}\\
	\hline
	\end{tabular}
\end{table}

In Table~\ref{tab:3aogp}, the groups $A_p(n,\theta)$ are known as \textit{Suzuki $p$-groups} of type A (see Example~\ref{exam:apn}); the group $P(\epsilon)$ was first given by Dornhoff~\cite{dornhoff1970imprimitive} (see Theorem~\ref{thm:2group}); and $q^{1+2m}_{+}$ is a generalization of $p^{1+2m}_{+}$, the \textit{extraspecial $p$-groups} of exponent $p$ with odd prime $p$.
We now explain how to generalize $p^{1+2m}_{+}$ to $q^{1+2m}_{+}$.

Let $q=p^n$ for an odd prime $p$.
For convenience, we denote a Sylow $p$-subgroup of a group $X$ by $X_p$.
It is known that $p^{1+2}_{+}\cong \SL(3,p)_p$ and $p^{1+2m}_{+}$ is the central product of $m$ copies of $p^{1+2}_{+}$.
A natural way to generalize such groups is defining $q^{1+2}_{+}$ to be $\SL(3,q)_p$ and defining $q^{1+2m}_{+}$ to be the central product of $m$ copies of $q^{1+2}_{+}$ (also see Definition~\ref{def:extra}).
A classical result of Winter~\cite{winter1972automorphism} says that $\Out(p^{1+2m}_{+})\cong\mathrm{CSp}(2m,p)\cong\Sp(2m,p){:}\bbZ_{p-1}$.
In Lemma~\ref{lem:nmpncgamma}, we prove that $\Aut(q^{1+2m}_{+})\cong p^{2mn^2}{:}\CGammaSp(2m,q)$, which shows some similarities between $p^{1+2m}_{+}$ and $q^{1+2m}_{+}$.
Thus, we name the group $q^{1+2m}_{+}$ the \textit{extraspecial $q$-group}.

For odd prime $p$, we can observe from Table~\ref{tab:3aogp} that a non-abelian $3$-orbit $p$-group is isomorphic to either $A_p(n,\theta)$ with $|\theta|=3$ or a quotient group of $q^{1+2m}_{+}$ for some $q=p^n$.
The following corollary is a criterion for deciding the quotient group $q^{1+2m}_+/U$ to be $3$-orbit.
\begin{corollary}\label{coro:cgammaspfinal}
	Let $N=q^{1+2m}_{+}$, and let $H=\Aut(N)^{\bfZ(N)}\cong \GammaL(1,q)$.
	Then $N/U$ is a $3$-orbit group with $U<\bfZ(N)$ if and only if $H_U$ acts transitively on the non-identity elements of $\bfZ(N)/U\cong\bfZ(N/U)$.
\end{corollary}

Moreover, we give some interesting properties of extraspecial $q$-groups.
We remark that an $\bbF_p$-subspace $U$ of $V=\bbF_{p^n}^+$ is called a \textit{subfield hyperplane} with respect to subfield $\bbF_{p^d}$ for $d\mid n$ if and only if $U$ is a hyperplane of $V$ as $\bbF_{p^d}$-spaces (also see Definition~\ref{def:subfield}). 

\begin{thm}\label{thm:isos}
	Let $q=p^n$ for odd prime $p$ and $n\geqs 1$.
	Then the following statements hold.
	\begin{enumerate}
		\item Let $U<\bfZ(q^{1+2m}_+)\cong\bbF_{q}^+$.
		Then $U$ is a subfield hyperplane of $\bfZ(q_{+}^{1+2m})\cong \bbF_{q}^+$ with respect to $\bbF_{q_0}$ if and only if 
		\[q^{1+2m}_{+}/U\cong (q_0)^{1+2mn/d}_{+},\]
		where $q_0=p^d$. 
		In particular, for any subgroup $U<\bfZ(q^{1+2m}_+)$ of order $p^{n-1}$, we have
		\[q^{1+2m}_{+}/U\cong p^{1+2mn}_{+}.\]
		\item $q^{1+2}_{+}\cong \SL(3,q)_p\cong \SU(3,q)_p$.
		\item $q^{1+2m}_{+}\cong \textbf{O}_p(G_\alpha)$, where $G_\alpha$ is the stabilizer of $\PSp(2m+2,q)$ acting naturally on $1$-dimensional subspaces of $\bbF_{q}^{2m+2}$.
	\end{enumerate}
\end{thm}

Let $N$ be a non-abelian $3$-orbit $p$-group with center $\bbZ_p$ and exponent $p$.
Theorems~\ref{thm:3orbits} and~\ref{thm:isos}\,(i) imply that $N$ is an extraspecial $p$-group $p^{1+2m}_+$ if $N$ is finite.
M{\"a}urer and Stroppel studied the structure of such infinite groups~\cite[Theorem 6.4]{maurer1997Groups}.

A transitive permutation group $G$ on a set $\Omega$ is said to be of \textit{rank} $r$ if $G$ has exactly $r$ orbits on $\Omega\times\Omega$.
The study of rank $3$ groups dates back to work of G.\,Higman~\cite{higman1964Finite} in 1964.
Note that the \textit{holomorph} $\Hol(N)=N{:}\Aut(N)$ of a group $N$ has rank $r$ if and only if $N$ is an $r$-orbit group.
Recently, Huang and the authors established a reduction theorem towards a classification of the finite rank $3$ permutation groups~\cite[Corollary 3]{huang2025finite}, of which case~(iii) comprises the groups with regular normal subgroups.
Hence, our classification of $3$-orbit groups provides the ``maximal'' groups in this case.

\begin{corollary}
	The holomorph $\Hol(N)\cong N{:}\Aut(N)$ is a finite rank $3$ permutation group if and only if $N$ is one of the groups in Table~\ref{tab:3aogp}.
\end{corollary}

The layout of the paper is as follows. 
In Section 2, we introduce some finite $3$-orbit groups, and then figure out the relations between representation theory and special $p$-groups.
The proof of Theorem~\ref{thm:transmod} will be given in Section~3.
Section~4 contains constructions and proofs of properties of groups in rows 1-7 in Table~\ref{tab:3aogp}.
We will discuss the relations between subfield hyperplanes and groups in the 8-th row in Table~\ref{tab:3aogp} in Section~5, and then we complete proofs of Theorems~\ref{thm:3orbits} and~\ref{thm:isos}.

We will introduce some known finite $3$-orbit groups in Section~\ref{sec:exam}, for example, $\{p,q\}$-groups, abelian groups and $2$-groups.
Then some basic properties of special $p$-groups of exponent $p$ will be given in Section~\ref{sec:special}.
A new family of $3$-orbit $p$-group, the \textit{extraspecial $q$-groups}, will be constructed in Section~\ref{sec:extra}.

The proof of Theorem~\ref{thm:transmod} will be given in Section~\ref{sec:mod}, which is a generalization of~\cite[Theorem 1.1]{li2024proof}.
With the aid of this theorem, we obtain the full automorphism group of the groups in rows 1-7 of Table~\ref{tab:3aogp} in Section~\ref{sec:oddp}.

We will discuss some properties of quotients of extraspecial $q$-groups in Section~\ref{sec:quoextra}, and then complete proofs of Theorems~\ref{thm:3orbits} and~\ref{thm:isos}.

\subsection*{Acknowledgments}
This work was supported by NNSFC grant no.~11931005.
The authors thank Eamonn~A. O'Brien for his advice on \Magma~computations, and are grateful to the anonymous referees for valuable comments and suggestions, which have helped improve the paper.

\section{$3$-orbit groups and special $p$-groups}\label{sec:pgroup}

In this section, we draw some preliminary properties and results on $3$-orbit groups.
As observed before, a finite $3$-orbit group is either a $p$-group of exponent at most $p^2$, or a $\{p,q\}$-group for two primes $p,q$.

\subsection{Some known finite $3$-orbit groups}\label{sec:exam}\

Finite $3$-orbit $\{p,q\}$-groups are classified in~\cite{devillers2014normal,dornhoff1970imprimitive,laffey1986Automorphism,maurer1997Groups}, and are described in the following Example~\ref{exam:pq}.
Recall that a \textit{primitive prime divisor} of $a^k-1$ is a prime $r$ such that $r\mid (a^k-1)$ but $r\nmid (a^i-1)$ for any $1\leqs i\leqs k-1$.

\begin{example}\label{exam:pq}
	Let $p$ and $q$ be two distinct primes such that $q$ is a primitive prime divisor of $p^{q-1}-1$.
	Set $X=\AGammaL(n,p^{q-1})=P{:}H$ with $P\cong\bbZ_p^{n(q-1)}$ and $H\cong\GammaL(n,p^{q-1})$.
	Then $H$ has a normal cyclic subgroup $L\cong\bbZ_{p^{q-1}-1}$ which is the center of $\GL(n,p^{q-1})$.
	Noting that $q\mid p^{q-1}-1$, let $Q$ be the unique cyclic subgroup of $L$ of order $q$.
	Then the $\{p,q\}$-group
	\[P{:}Q\cong\bbZ_{p}^{n(q-1)}{:}\bbZ_q\]
	is a normal subgroup of $X=P{:}H=\AGammaL(n,p^{q-1})$.
\end{example}

Let $N=P{:}Q$ as defined in Example~\ref{exam:pq}.
We make the following observations on $N$.
\begin{enumerate}
	\item Viewing $P\cong\bbZ_p^{n(q-1)}$ as a vector space, $Q$ acts on $P$ via scalar multiplications by definition, and so $P$ is a homogeneous $\bbF_pQ$-module.
	Thus, $N$ is a homogeneous Frobenius group, see~\cite[Theorem 2]{laffey1986Automorphism}.
	Then, by~\cite[Theorem 1.1]{wang2020automorphism}, we can easily obtain that $\Aut(N)\cong X=\AGammaL(n,p^{q-1})$.

	\item Elements in $N$ of order $p$ are conjugate under $\GammaL(n,p^{q-1})$.
	Sylow $q$-subgroups of $N$ are conjugate.
	The field automorphism $\phi$ of $\bbF_{p^{q-1}}$ acts on $Q$ faithfully, and so non-identity elements of $Q$ are conjugate under $\langle \phi\rangle$.
	Thus, elements in $N$ of order $q$ are conjugate under $\Aut(N)$, and $N=P{:}Q$ is a $3$-orbit group.
\end{enumerate}

By~\cite{devillers2014normal,dornhoff1970imprimitive,laffey1986Automorphism,maurer1997Groups}, each $3$-orbit group of order divisible by primes $p,q$ is isomorphic to a homogeneous Frobenius group $\bbZ_{p}^{n(q-1)}{:}\bbZ_q$, as described in Example~\ref{exam:pq}.
\begin{lemma}\label{lem:pq}
	A finite $3$-orbit group $N$ is not a $p$-group if and only if $N\cong\bbZ_{p}^{n(q-1)}{:}\bbZ_q$, as described in Example~$\ref{exam:pq}$.
	In this case, $\Aut(N)\cong\AGammaL(n,p^{q-1})$.
\end{lemma}

Obviously, each homocyclic group $\bbZ_{p^2}^n$ is a $3$-orbit group, and
\[\Aut(\bbZ_{p^2}^n)=\GL(n,\bbZ/p^2\bbZ)\cong p^{n^2}.\GL(n,p).\]
Conversely, every finite abelian $3$-orbit group is homocyclic of exponent $p^2$.
\begin{lemma}\label{lem:ab}
	A finite abelian group $N$ is a $3$-orbit group if and only if $N\cong\bbZ_{p^2}^n$ for some prime $p$ and integer $n$.
\end{lemma}

We note that Schwachh{\"o}fer and Stroppel~\cite{schwachhofer1999Finding} studied the structure of $k$-orbit abelian groups (especially for infinite groups) for finite $k$.

A non-abelian $2$-group with more than one involution is called a \textit{Suzuki $2$-group} if it has a cyclic group of automorphisms acting transitively on involutions.
D.\,Higman classified all Suzuki $2$-groups and divided them into four types A-D in~\cite{higman1963Suzuki}.
The definition of each type of Suzuki $2$-groups can be extended to odd primes, which are called \textit{Suzuki $p$-groups}.
The following is the definition of Suzuki $p$-group of type A, see~\cite[Section 46]{berkovich2008Groupsa}.
\begin{example}\label{exam:apn}
	Let $\theta\in\Gal(\bbF_{p^n}/\bbF_{p})$.
	For $a,b\in\bbF_{p^n}$, define
	\[A_p(n,\theta)=\langle M(a,b,\theta) \mid a,b\in\bbF_{p^n}\rangle,\]
	where $M(a,b,\theta)$ is a matrix of the form
	\[M(a,b,\theta)=\begin{pmatrix}1&a&b\\0&1&a^\theta\\0&0&1\end{pmatrix}.\]
\end{example}

We remark that groups $A_2(n,\theta)$ with $|\theta|\neq 1$ of odd order form the class of Suzuki $2$-groups of type A given by D.\,Higman, which are $3$-orbit groups.
The following proposition can be easily verified, where the first part was given in~\cite[Lemma~4.2]{li2024proof}.

\begin{proposition}\label{prop:antheta}
	Let $N=A_p(n,\theta)$.
	Then the following statements hold.
	\begin{enumerate}
		\item When $p=2$ and $|\theta|>1$ is odd, $N$ is a $3$-orbit $2$-group with $\Aut(N)\cong 2^{n^2}{:}\GammaL(1,2^n)$.
		\item Let $\lambda$ be a generator of $\bbF_{p^n}$, and let $D=\mathrm{diag}(\lambda^{-1},1,\lambda^\theta)$.
		Then $D$ induces an automorphism $\xi$ of $N$ such that
		\[\xi:M(a,b,\theta)\mapsto D^{-1}M(a,b,\theta)D=M(a\lambda,b\lambda\lambda^\theta,\theta).\]
	\end{enumerate}
\end{proposition}

We remark that $A_2(2f+1,\theta)$ for $x^\theta=x^{2^{f+1}}$ is isomorphic to a Sylow $2$-subgroup of Suzuki group $\Sz(2^{2f+1})$, see~\cite[Proposition 13.4\,(vi)]{carter1972Simple}.
Hence, Sylow $2$-subgroups of $\Sz(2^{2f+1})$ are $3$-orbit $2$-groups.

Finite $3$-orbit $2$-groups are classified independently by Bors and Glasby~\cite{bors2024Finite}, and the authors~\cite{li2024proof}.
We record the results below.
\begin{theorem}\label{thm:2group}
	Suppose that $N$ is a non-abelian $3$-orbit $2$-group.
	Then one of the following statement holds.
	\begin{enumerate}
		\item $N\cong A_2(n,\theta)$ for $|\theta|\neq 1$ of odd order, and $\Aut(N)\cong 2^{n^2}{:}\GammaL(1,2^n)$.
		\item $N\cong \SU(3,2^n)_2$ and $\Aut(N)\cong 2^{2n^2}{:}\GammaL(1,2^{2n})$.
		\item $N\cong P(\epsilon)$ and $\Aut(N)\cong 2^{18}{:}(7{:}9)$, where
		\[\begin{aligned}
			P(\epsilon)=&\langle x_i,z_j\mbox{ for $1\leqs i\leqs 6$ and $1\leqs j\leqs 3$ }|z_j^2=[x_i,z_j]=[z_k,z_{\ell}]=1\mbox{ for each $i,j,k,\ell$,}\\
			&x_1^2=x_3^2=z_2,\ x_2^2=z_2z_3,\ x_4^2=z_3,\ x_5^2=z_1z_2z_3,\ x_6^2=z_3,\\
			&[x_1,x_2]=[x_3,x_5]=[x_3,x_6]=z_1z_2,\ [x_1,x_3]=z_1z_3,\ [x_1,x_4]=z_3,\\
			&[x_1,x_5]=[x_3,x_4]=[x_5,x_6]=z_2,\ [x_1,x_6]=1,\ [x_2,x_6]=z_1z_2z_3,\\
			&[x_2,x_3]=[x_2,x_4]=[x_4,x_6]=z_1,\ [x_2,x_5]=[x_4,x_5]=z_2z_3.\rangle
		\end{aligned}\]
	\end{enumerate}
\end{theorem}

We remark that groups in parts~(ii) and (iii) are Suzuki $2$-groups of type B in Higman's notation, see~\cite[page 706]{dornhoff1970imprimitive} and~\cite[Theorem 1.2]{bors2024Finite}.

\subsection{Special $p$-groups of exponent $p$}\label{sec:special}\

Let $N$ be a non-abelian $3$-orbit $p$-group with $p$ an odd prime.
Then non-identity elements in the center $\bfZ(N)$ form an automorphism orbit, and elements in $N\setminus\bfZ(N)$ form a different automorphism orbit.
It follows that both $\bfZ(N)$ and $N/\bfZ(N)$ are elementary abelian groups, and $\bfZ(N)=\Phi(N)=N'$.
So $N$ is a special $p$-group of exponent either $p$ or $p^2$.
If the exponent of $N$ is $p^2$, then all elements of $N$ of order $p$ form an automorphism orbit, and so $N$ is abelian by~\cite[Corollary 3]{shult1969finite}.
This leads to the following lemma.
\begin{lemma}\label{lem:3aop}
	Let $N$ be a non-abelian $3$-orbit $p$-group with $p$ an odd prime.
	Then $N$ is a special $p$-group with exponent $p$.
\end{lemma}

For a non-abelian special $p$-group $N$, $\bfZ(N)=N'=\Phi(N)$ and $N/\Phi(N)$ are elementary abelian $p$-groups, and thus they can be viewed as vector spaces over field $\bbF_p$.
Let $A=\Aut(N)$.
Then the actions of $A$ on $\Phi(N)$ and $N/\Phi(N)$ are linear actions.
It yields that $\Phi(N)$ and $N/\Phi(N)$ can be viewed as $\bbF_pA$-modules.
For elements $x,y,z\in N$, we have $[x,yz]=[x,y]^z[x,z]=[x,y][x,z]$ as $\bfZ(N)=N'$.
Thus, the commutator operator $[\cdot,\cdot]$ on $N$ resembles an alternating bilinear form.
Let $V$ be an $\F G$-module.
The \textit{exterior square} $\Lambda_\F^2(V)$ is the quotient $\F G$-module $(V\otimes V)/S$, where $S$ is spanned by vectors $v\otimes v$ for all $v\in V$.
The exterior square has a universal property associated with alternating bilinear form.
The following lemma is given in~\cite[Lemma 3.1]{li2024proof}.

\begin{lemma}\label{lem:spext}
	Let $N$ be a non-abelian special $p$-group with $p$ prime such that $N/\bfZ(N)\cong \bbZ_p^n$.
	Set $V=N/\bfZ(N)$ and $M=\bfZ(N)$. 
	Let $K$ and $L$ be the kernels of $\Aut(N)$ acting on $V$ and $M$, respectively.
	Define $G=A/K$.
	Then the following statements hold.
	\begin{enumerate}
		\item $K$ is an elementary abelian $p$-group of order $|\bfZ(N)|^n$.
		\item $K\leqs L$, and then both $V$ and $M$ are naturally $\bbF_p G$-modules.
		\item $N$ is an $r$-orbit group with $r=o(V)+o(M)-1$, where $o(V)$ and $o(M)$ are numbers of $G$-orbits on $V$ and $M$, respectively.
		\item $M\cong\Lambda^2_{\bbF_p}(V)/W$ for some $\bbF_p G$-submodule $W$ of $\Lambda^2_{\bbF_p}(V)$.
	\end{enumerate}
\end{lemma}

Recall that, for an odd prime $p$, a non-abelian $3$-orbit $p$-group $N$ is a special $p$-group of exponent $p$.
Let $\calH_{n,p}$ denote the $n$-generator free group in the variety of groups that have exponent $p$ and nilpotency class $2$.
Then $N$ is a quotient of $\calH_{n,p}$ when $N/\Phi(N)\cong \bbZ_p^n$.
We remark that $\calH_{n,p}$ is the so-called the \textit{universal generalized Heisenberg group} in~\cite{maurer1997Groups}, and it is a finite special $p$-group with $|\bfZ(\calH_{n,p})|=p^{\binom{n}{2}}$ (see~\cite[Page 107]{berkovich2011Groups}).
We formulate an equivalent definition below.
\begin{definition}\label{def:hei}
	Let $p$ be an odd prime, and let $n$ be a positive integer.
	Define
	\[\calH_{n,p}=\langle \alpha_1,...,\alpha_n|\alpha_i^p=[\alpha_i,[\alpha_j,\alpha_k]]=1\mbox{ for $i,j,k=1,...,n$}.\rangle\]
\end{definition}

Let $\calH=\calH_{n,p}$ and let $W<\bfZ(\calH)$.
The automorphism group of $\calH/W$ for $W<\bfZ(\calH)$ is described in \cite[Lemma~3]{glasby2011groups} which is based on~\cite[Theorem~2.10]{obrien1990group}.
We state the results as below.
\begin{lemma}\label{lem:oddp}
	Suppose that $N=\calH/W$ where $\calH=\calH_{n,p}$ and $W<\bfZ(\calH)$.
	Let $V=\calH/\bfZ(\calH)=N/\bfZ(N)$, and let $M=\bfZ(N)=\bfZ(\calH)/W$.
	Define $X=\Aut(\calH)^V$ and $G=\Aut(N)^V$.
	Then the following statements hold.
	\begin{enumerate}
		\item $X=\Aut(V)=\GL(V)\cong\GL(n,p)$, and $\bfZ(\calH)\cong \Lambda^2_{\bbF_p}(V)$ as an $\bbF_pX$-module.
		\item $G\leqs X$ is equal to $X_{W}$, the maximal subgroup of $X$ that stabilizes the subspace $W$ of $\bfZ(\calH)$.
	\end{enumerate}
\end{lemma}

Using the notations above, Lemma~\ref{lem:spext}~(i) implies that
\[\Aut(\calH)=K.X\cong p^{\frac{n^2(n-1)}{2}}.\GL(V),\]
where $K\cong p^{\frac{n^2(n-1)}{2}}$ is the kernel of $X$ acting on $\calH/\bfZ(\calH)$.
The following simple lemma shows that the group extension $K.X$ is split.
\begin{lemma}\label{lem:split}
	Suppose that $N=\calH/W$ where $\calH=\calH_{n,p}$ and $W<\bfZ(\calH)$.
	Let $V=N/\bfZ(N)$, and let $K$ be the kernel of $\Aut(N)$ acting on $V$.
	Then the group extension $\Aut(N)=K.\Aut(N)^{V}$ is split.
\end{lemma}
\begin{proof}
	First, we consider $N=\calH$.
	Then $K\cong p^{\frac{n^2(n-1)}{2}}$ and $\Aut(\calH)^{V}\cong\GL(V)\cong\GL(n,p)$ by Lemma~\ref{lem:oddp}.
	Recall that $\calH=\langle\alpha_1,...,\alpha_n\rangle$.
	Let $k\in\bbF_p^\times$ be a primitive element.
	Then there exists $\psi\in\Aut(\calH)$ of order $p-1$ such that 
	\[\psi:\alpha_i\mapsto x_i\alpha_i^k \mbox{, for some $x_i\in\bfZ(\calH)$ with $i=1,...,n$.}\]
	Since $K$ is the kernel of $\Aut(\calH)$ acting on $V=\calH/\bfZ(\calH)$, for each $\varphi\in K$ there exist $a_1,...,a_n\in \bfZ(\calH)$ such that 
	\[\varphi: \alpha_i\mapsto a_i\alpha_i,\mbox{ for }i=1,...,n.\]
	Note that $a_i$ is a product of some $[\alpha_i,\alpha_j]$'s and 
	\[[\alpha_i,\alpha_j]^\psi=[\alpha_i^\psi,\alpha_j^\psi]=[\alpha_i,\alpha_j]^{k^2}.\]
	Thus, we have $a_i^\psi=a_i^{k^2}$ for each $i,...,n$.
	Let $s\in\bbF_p^\times$ such that $sk=1$.
	Then $\alpha_i^{\psi^{-1}}=x_i^{-s^3}\alpha_i^s$ and
	\[\alpha_i^{\psi^{-1}\varphi\psi}=(x_i^{-s^3}\alpha_i^s)^{\varphi\psi}=(x_i^{-s^3}a_i^s\alpha_i^s)^{\psi}=(a_i^s)^\psi \alpha_i=a_i^k\alpha_i.\]
	It follows that $\varphi^\psi=\varphi^k$, and hence $C_{K}(\psi)=1$.

	Let $\phi\in \Aut(\calH)$.
	Then $[\phi K,\psi K]=1$ since $\langle \psi K\rangle=\bfZ(\Aut(\calH)/K)\cong \bbZ(\GL(V))$.
	Note that $k\neq 1$, there exists $t\in\bbF_p^\times$ such that $t(1-k)=1$.
	Let $\varphi=\left[\psi^{\phi^{-1}},\phi\right]^{t}\in K$, and let $\phi'=\varphi\phi$.
	Then 
	\[\begin{aligned}
		[\phi',\psi]&=[\varphi,\psi]^\phi[\phi,\psi]=\left((\varphi^{-1})^{\psi}\varphi\right)^\phi[\phi,\psi]=\left(\varphi^{1-k}[\phi,\psi]^{\phi^{-1}}\right)^\phi\\
		&=\left(\left[\psi^{\phi^{-1}},\phi\right]\left[\phi,\psi^{\phi^{-1}}\right]\right)^\phi=1.
	\end{aligned}\]
	Thus, for any $\phi\in\Aut(\calH)$, there exists $\phi'\in\Aut(\calH)$ such that $\phi'K=\phi K$ and $[\phi',\psi]=1$.
	It follows that 
	\[|C_{\Aut(\calH)}(\psi)|\geqs|\Aut(\calH)^V|=|\GL(V)|.\]
	Since $C_{K}(\psi)=1$ and $\Aut(\calH)\cong K.\GL(V)$, we have $C_{\Aut(\calH)}(\psi)\cong\GL(V)$ is a complement of $K$ in $\Aut(\calH)$.

	Now we consider the general case that $N=\calH/W$.
	Let $X$ be a complement of $K$ in $\Aut(\calH)$.
	Then $X_W$ admits a group of automorphisms of $N$, and is a complement of $K$ in $\Aut(N)$ where $K$ is the kernel of $\Aut(N)$ acting on $V$.
	Therefore, the group extension $\Aut(N)=K.\Aut(N)^V$ is split as $\Aut(N)^V\cong X_W$ by Lemma~\ref{lem:oddp}\,(ii).
\end{proof}

Every $g\in\GL(n,p)$ induces an automorphism $\sigma\in\Aut(\calH)$ by Lemma~\ref{lem:oddp}\,(i), where $\calH=\calH_{n,p}$.
Thus, we have a straightforward isomorphism $\calH/W\cong \calH/W^\sigma$ for any $W<\bfZ(\calH)$, which yields the following lemma.
\begin{lemma}\label{lem:oddpiso}
	Let $\calH=\calH_{n,p}$, and let $W_1,W_2$ be two subgroups of $\bfZ(\calH)$.
	View $W_1$ and $W_2$ as two subspaces of $\Lambda^2_{\bbF_p}(V)$ where $V=\calH/\bfZ(\calH)$.
	Then $\calH/W_1\cong \calH/W_2$ if there exists $g\in \GL(V)$ such that $W_1^g=W_2$ in $\Lambda^2_{\bbF_p}(V)$.
\end{lemma}

\subsection{Extraspecial $q$-groups}\label{sec:extra}\

Let $p$ be an odd prime.
We introduce a new family of special $p$-groups in this part.
Recall that the \textit{extraspecial $p$-group} of exponent $p$, denoted by $p^{1+2m}_{+}$, is the central product of $m$ copies of $\SL(3,p)_p$.
This leads us to start with the Sylow $p$-subgroups of $\SL(3,p^n)$.
\begin{definition}
	Let $q=p^n$ for odd prime $p$.
	The \textit{Heisenberg group} (over finite field $\bbF_q$), denoted by $q^{1+2}_{+}$, is a subgroup of $\SL(3,q)$ consisting of matrices of form $\begin{pmatrix}1&a&c\\0&1&b\\0&0&1\end{pmatrix}$ for all $a,b,c\in\bbF_q$.
\end{definition}

We remark that the Heisenberg group is a non-abelian special $p$-group.
When $p=q$, it is exactly the extraspecial $p$-group of plus type, and our notation fits $p^{1+2}_{+}$.
\begin{lemma}\label{lem:q12}
	Let $N=q^{1+2}_{+}$ for $q=p^n$ where $p$ is an odd prime.
	Then $p^{2n^2}{:}\GammaL(2,q)\lesssim\Aut(N)$ and $N$ is a $3$-orbit group.
\end{lemma}
\begin{proof}
	For $x,y,z\in\bbF_q$, let
	\[A(x)=\begin{pmatrix}1&x&0\\0&1&0\\0&0&1\end{pmatrix}, B(y)=\begin{pmatrix}1&0&0\\0&1&y\\0&0&1\end{pmatrix}\mbox{ and }C(z)=\begin{pmatrix}1&0&z\\0&1&0\\0&0&1\end{pmatrix}.\]
	Then
	\begin{enumerate}
		\item[(a)] $[A(x),B(y)]=C(xy)$, $[A(x),A(y)]=[B(x),B(y)]=1$; and 
		\item[(b)] $[A(x_1),B(y_1)][A(x_2),B(y_2)]=[A(x_1+x_2),B(y_1+y_2)]$ for $x_1,x_2,y_1,y_2\in\bbF_q$.
	\end{enumerate}
	
	Since $N$ is a special $p$-group of exponent $p$, it is isomorphic to $\calH/W$ with $\calH=\calH_{2n,p}$ and $W<\bfZ(\calH)$.
	Let $\alpha(x)$ and $\beta(y)$ be preimages of $A(x)$ and $B(y)$ in $\calH$ for $x,y\in\bbF_q$, respectively.
	Note that each element in $N$ can be expressed as a product of elements $A(x)$, $B(y)$ and $C(z)$ for some $x,y,z\in\bbF_q$.
	With the equalities in (a), we have that $C(z)\in\langle  A(x),B(y):x,y\in \bbF_q\rangle$.
	This yields that 
	\[N=\langle A(x),B(y):x,y\in \bbF_q\rangle\mbox{ and }\calH=\langle\alpha(x),\beta(y):x,y\in\bbF_q\rangle.\]
	Let $X=\Aut(\calH)^{V}$ with $V=\calH/\bfZ(\calH)$.
	By Lemma~\ref{lem:oddp}\,(i), $X= \GL(V)\cong\GL(2n,p)$.
	Note that $V$ can be viewed as a $2$-dimensional $\bbF_q$-space $V_q$ with basis $\alpha(1)\bfZ(\calH)$ and $\beta(1)\bfZ(\calH)$.
	Let $G\cong \GammaL(2,p^n)$ be the subgroup of $X$ acting naturally on $V_q$.
	
	Now we prove that $W^g=W$ for $g\in G$, and then $G\lesssim \Aut(N)^V$ by Lemma~\ref{lem:oddp}\,(ii).
	Note that each element in $\bfZ(\calH)=\calH'$ is a product of elements of forms $[\alpha(x),\alpha(y)]$, $[\beta(x),\beta(y)]$ and $[\alpha(x),\beta(y)]$ for $x,y\in\bbF_q$.
	\begin{enumerate}
		\item Since $[A(x),A(y)]=1$, we have $[\alpha(x),\alpha(y)]\in W$.
		Assume that $(\alpha(1)\bfZ(\calH))^g=\alpha(a)\beta(b)\bfZ(\calH)$ for some $a,b\in\bbF_q$.
		Then
		\[[\alpha(x),\alpha(y)]^g=[\alpha(ax)\beta(bx),\alpha(ay)\beta(by)].\]
		We obtain that $[\alpha(x),\alpha(y)]^g\in W$ since
		\[\begin{aligned}
			[A(ax)B(bx),A(ay)B(by)]&=[A(ax),B(by)][B(bx),A(ay)]\\
			&=C(ax\cdot by-bx\cdot ay)=1.
		\end{aligned}\]
		\item Similarly, $[\beta(x),\beta(y)]\in W$ and $[\beta(x),\beta(y)]^g\in W$ for any $x,y\in\bbF_q$ and $g\in G$.
		\item We assume that $w\in W$ with $w=\prod_{i=1}^k[\alpha(x_i),\beta(y_i)]$ for some $x_i,y_i\in\bbF_q$.
		Notice that $\sum_{i=1}^kx_iy_i=0$ since $\prod_{i=1}^k[A(x_i),B(y_i)]=C\left(\sum_{i=1}^nx_iy_i\right)$.
		Let $a,b,c,d\in\bbF_q$ such that
		\[(\alpha(1)\bfZ(\calH))^g=\alpha(a)\beta(b)\bfZ(\calH)\mbox{ and }(\beta(1)\bfZ(\calH))^g=\alpha(c)\beta(d)\bfZ(\calH).\]
		By calculations, we have $w^g$ equals to
		\[\prod_{i=1}^k[\alpha(x_i),\beta(y_i)]^g=\prod_{i=1}^k[\alpha(x_ia)\beta(x_ib),\alpha(y_ic)\beta(y_id)]=\prod_{i=1}^k s_i \cdot \prod_{i=1}^k t_i,\]
		where $s_i=[\alpha(x_ia),\alpha(y_ic)][\beta(x_ib),\beta(y_id)]$, and $t_i= [\alpha(x_ia),\beta(y_id)][\beta(x_ib),\alpha(y_ic)]$.
		By (i) and (ii), each $s_i$ lies in $W$.
		The image of $\prod_{i=1}^k t_i$ in $N$ is equal to
		\[\prod_{i=1}^k[A(x_ia),B(y_id)][B(x_ib),A(y_ic)]=C\left((ad-bc)\cdot \sum_{i=1}^kx_iy_i\right).\]
		Since $\sum_{i=1}^kx_iy_i=0$, the image of $\prod_{i=1}^kt_i$ in $N$ is trivial.
		Hence, the image of $w^g$ in $N$ is trivial, that is, $w^g\in W$.
	\end{enumerate}

	Hence, $W^G=W$ and $G\lesssim \Aut(N)^V$ by Lemma~\ref{lem:oddp}\,(ii).
	Thus, $\Aut(N)\geqs K.G\cong p^{2n^2}.\GammaL(2,q)$ where $K$ is the kernel of $\Aut(N)$ acting on $N/\bfZ(N)$ by Lemma~\ref{lem:spext}\,(i).
	In particular, the extension $K.G$ is split by Lemma~\ref{lem:split}, and then $p^{2n^2}{:}\GammaL(2,q)\lesssim \Aut(N)$.

	Note that there exists an automorphism $\xi$ of $N$ admitted by conjugating of diagonal matrix $\mathrm{diag}(\lambda,1,1)$ where $\lambda$ is a generator of $\bbF_{q}$.
	By definition we observe that $\langle\xi\rangle$ acts regularly on non-identity elements of $\bfZ(N)$.
	In addition, $G\cong \GammaL(2,q)$ acts transitive on non-zero vectors of $V=N/\bfZ(N)\cong\bbF_p^{2n}$.
	Hence $\Aut(N)$ acts transitively on non-identity elements of both $N/\bfZ(N)$ and $\bfZ(N)$.
	By Lemma~\ref{lem:spext}\,(iii), $N$ is a $3$-orbit group.
\end{proof}

Note that $\GammaL(2,q)$ contains the Singer cyclic which is a cyclic group of order $q^2-1$ acting regularly on non-identity elements of $\bbF_q^2$.
This yields a solvable subgroup of $\Aut(N)$ which acts on $N$ with $3$ orbits.
Immediately, we obtain the following lemma, which is Theorem~\ref{thm:isos}\,(ii).
\begin{lemma}\label{lem:susl}
	Let $q=p^n$ with odd prime $p$.
	Then \[q^{1+2}_{+}\cong\SL(3,q)_p\cong \SU(3,q)_p.\]
\end{lemma}
\begin{proof}
	Let $N=q^{1+2}_{+}$, and let $A\cong p^{2n^2}{:}\GammaL(2,q)$ be the subgroup of $\Aut(N)$ given in Lemma~\ref{lem:q12}.
	Then there exists $\varphi\in A$ of order $q^2-1$ such that $\varphi$ induces a transitive action on non-identity elements of $N/\bfZ(N)$.
	Let $\xi\in\Aut(N)$ induced by conjugating the diagonal matrix $\mathrm{diag}(\mu,1,1)$ for some generator $\mu$ of $\bbF_q^\times$.
	We remark that $\xi\in A$ and $\xi$ acts transitively on non-identity elements of $\bfZ(N)$.
	Let $K\cong p^{2n^2}$ be the kernel of $A$ acting on $N/\bfZ(N)$.
	Then $\langle K,\varphi,\xi\rangle$ is a solvable subgroup of $A\leqs\Aut(N)$, and it acts on elements of $N$ with exactly $3$ orbits by Lemma~\ref{lem:spext}\,(iii).
	Thus, $\Aut(N)$ has a solvable subgroup with exactly $3$ orbits on $N$.

	Let $H$ be a Sylow $p$-subgroup of $\SU(3,q)$ generated by matrices $M(a,b)\in\GL(3,q^2)$ such that $b+b^q+a^{1+q}=0$ and 
	\[M(a,b)= \begin{pmatrix}1&a&b\\0&1&a^q\\0&0&1\end{pmatrix},\mbox{ see~\cite[page 249]{dixon1996Permutation}}.\]
	Let $D=\mathrm{diag}(\lambda^{-q},\lambda^{1-q},\lambda)$, where $\langle\lambda\rangle=\bbF_{q^2}^\times$.
	Then $D^{-1}M(a,b)D=M(q\lambda,b\lambda^{q+1})$ and
	\[b\lambda^{q+1}+(b\lambda^{q+1})^q+(a\lambda)(a\lambda)^q=\lambda^{1+q}(b+b^q+a^{1+q}).\]
	Hence, $D$ induces an automorphism $\xi'\in\Aut(H)$ such that $\langle\xi'\rangle$ acts transitively on non-identity elements of both $H/\bfZ(H)$ and $\bfZ(H)$.
	By Lemma~\ref{lem:spext}\,(i) and (iii), $\Aut(H)$ has a solvable subgroup with exactly $3$ orbits on $H$.

	Finite groups which have solvable groups of automorphisms with exactly $3$ orbits have been classified in~\cite{dornhoff1970imprimitive}.
	In particular, there exists a unique such group (up to isomorphism) of order $p^{3n}$ with center of order $p^n$.
	Note that both $q^{1+2}_{+}$ and $\SU(3,q)_p$ has order $p^{2n}$ with center of order $p^n$.
	Therefore, we have $q^{1+2}_{+}\cong\SL(3,q)_p\cong \SU(3,q)_p$.
\end{proof}

With Lemma~\ref{lem:susl}, we can reformulate Dornhoff's Theorem given in~\cite{dornhoff1970imprimitive}.
\begin{theorem}[Dornhoff]\label{thm:revdorn}
	Let $N$ be a finite group.
	If there exists a solvable subgroup of $\Aut(N)$ acting on $N$ with exactly $3$ orbits, then one of the following statement holds.
	\begin{enumerate}
		\item $N$ is an abelian group isomorphic to $\bbZ_p^n$ or $\bbZ_{p^2}^n$ for prime $p$.
		\item $N$ is a $\{p,q\}$-group given in Example~\ref{exam:pq}.
		\item $N$ is a non-abelian special $2$-group isomorphic to one of $A_2(n,\theta)$ (with $|\theta|\neq 1$ odd), $\SU(3,2^n)_2$ and $P(\epsilon)$.
		\item $N\cong q^{1+2}_{+}\cong\SU(3,q)_p\cong\SL(3,q)_p$ for odd prime power $q=p^n$.
		\item $N\cong 3^{1+4}_{+}$.
	\end{enumerate}
\end{theorem}

Now we extend the definition of extraspecial $p$-groups.
\begin{definition}\label{def:extra}
	For $q=p^n$ with odd prime $p$.
	The \textit{extraspecial $q$-group} is defined by
	\[q^{1+2m}_{+}=q^{1+2}_{+}\circ\cdots\circ q^{1+2}_{+}=(q^{1+2}_{+})^m/C,\]
	where $C$ consists of elements $(X_1,...,X_m)\in \bfZ(q^{1+2}_{+})^m$ such that $X_1\cdots X_m=1$.
\end{definition}

It is straightforward to see that $q^{1+2m}_{+}$ is a non-abelian special $p$-group of exponent $p$ such that $\bfZ(q^{1+2m}_{+})\cong\bbZ_p^n$ and $q^{1+2m}_{+}/\bfZ(q^{1+2m}_{+})\cong\bbZ_p^{2mn}$.
Let $\pi_i$ be the natural homomorphism from $q^{1+2}_{+}$ to the image of the $i$-th component in $q^{1+2m}_{+}$ for $i=1,...,m$.
Remark that each matrix $Z\in \bfZ(q^{1+2}_{+})$ satisfies $\pi_i(Z)=\pi_j(Z)$ for $i,j=1,...,m$.
Then we define the following notations for $i=1,...,m$ and $a,b,c\in\bbF_q$:
\[
	E_i(a)=\pi_i\begin{pmatrix}1&a&0\\0&1&0\\0&0&1\end{pmatrix}\mbox{, }
	F_i(b)=\pi_i\begin{pmatrix}1&0&0\\0&1&b\\0&0&1\end{pmatrix}\mbox{, and }
	Z(c)=\pi_i\begin{pmatrix}1&0&c\\0&1&0\\0&0&1\end{pmatrix}.
\]
The following equalities can be verified by simple calculations on matrices.
\begin{proposition}\label{prop:equalnmq}
	Using notations defined above.
	\begin{enumerate}
		\item $\pi_i\begin{pmatrix}1&a&c\\0&1&b\\0&0&1\end{pmatrix}=Z(c)F_i(b)E_i(a)$
		\item $E_i(a)E_i(b)=E_i(a+b),\ F_i(a)F_i(b)=F_i(a+b)\mbox{ and }\ Z(a)Z(b)=Z(a+b)$;
		\item $[E_i(a),E_j(b)]=[E_i(a),F_j(b)]=[F_i(a),F_j(b)]=1$ for $i\neq j$;
		\item $[E_i(a),F_i(b)]=Z(ab)$;
		\item $E_i(a)F_i(b)=F_i(b)E_i(a)[E_i(a),F_i(b)]=Z(ab)F_i(b)E_i(a)$.
	\end{enumerate}
\end{proposition}

With above properties, we obtain the following lemma about $\Aut(q^{1+2m}_{+})$.
\begin{lemma}\label{lem:autnmq}
	Let $N=q^{1+2m}_{+}$ with $q=p^n$ for odd prime $p$.
	Then $\Aut(N)$ has a subgroup $A$ such that $A^V\cong \CGammaSp(2m,q)$ and $A^{\bfZ(N)}\cong\GammaL(1,q)$, where $V=N/\bfZ(N)$.
\end{lemma}
\begin{proof}
	Since $N$ is a non-abelian special $p$-group of exponent $p$, $N=\calH/W$ with $\calH=\calH_{2mn,p}$ and $W\leqs\Phi(\calH)$.
	Identify $V=N/\bfZ(N)$ with $\calH/\bfZ(\calH)$, and let $X=\Aut(\calH)^{V}=\GL(V)$ as in Lemma~\ref{lem:oddp}.
	For convenience, let $\overline{\alpha}=\alpha\bfZ(\calH)\in V$ for $\alpha\in \calH$.
	Let $e_i(a)$ and $f_j(b)$ be preimages of $E_i(a)$ and $F_i(b)$ in $\calH$, respectively.
	Hence, $V$ can be viewed as an $\bbF_q$-symplectic space $V_q$ equipped with the symplectic form $\calB$ and a basis
	\[\overline{e_1(1)},...,\ \overline{e_m(1)},\ \overline{f_1(1)},...,\ \overline{f_m(1)}\]
	such that $\overline{e_i(1)}$ and $\overline{f_i(1)}$ form a hyperbolic pair for each $i$.
	Let $G$ be the subgroup of $X$ isomorphic to $\CGammaSp(2m,q)$ acting naturally as a semilinear conformal symplectic group on $V_q$.
	We claim that $w^g\in W$ for $g\in G$ and $w\in W$, then $G\lesssim\Aut(N)^V$ by Lemma~\ref{lem:oddp}\,(ii).
	Note that there exist some $\alpha_i,\beta_i\in \calH$ such that $w=\prod_{i=1}^k[\alpha_i,\beta_i]$.
	We prove the claim step by step.

	\textbf{Step 1}: we show that the image of $[\alpha,\beta]$ in $N$ is equal to $Z\left(\calB\left(\overline{\alpha},\overline{\beta}\right)\right)$ for any $\alpha,\beta\in\calH$.
	Assume that
	\[\overline{\alpha}=\prod_{i=1}^m\overline{e_i(a_i)f_i(b_i)}\mbox{ and }\overline{\beta}=\prod_{i=1}^m\overline{e_i(c_i)f_i(d_i)},\]
	where $a_i,b_i,c_i,d_i\in\bbF_q$.
	Then the image of $[\alpha,\beta]$ in $N$ equals
	\[\begin{aligned}
		\left[\prod_{i=1}^m(E_i(a_i)F_i(b_i)),\prod_{i=1}^m(E_i(c_i)F_i(d_i))\right]&=\prod_{i=1}^m[E_i(a_i),F_i(d_i)][F_i(b_i),E_i(c_i)]\\
		&=\prod_{i=1}^mZ(a_id_i-b_ic_i)=Z\left(\sum_{i=1}^m(a_id_i-b_ic_i)\right)\\
		&=Z\left(\calB\left(\overline{\alpha},\overline{\beta}\right)\right).
	\end{aligned}
	\]

	\textbf{Step 2}: we claim that $w\in W$ if and only if $\sum_{i=1}^k\calB\left(\overline{\alpha_i},\overline{\beta_i}\right)=0$.
	By the previous step, we obtain that the image of $w$ in $N$ is equal to $Z\left(\sum_{i=1}^k\calB\left(\overline{\alpha_i},\overline{\beta_i}\right)\right)$.
	This yields the claim since $w\in W$ if and only if $w$ has a trivial image in $N$.

	\textbf{Step 3}:
	we are now ready to show that $w^g\in W$.
	Note that $w^g=\prod_{i=1}^k[\alpha_i^{\widetilde{g}},\beta_i^{\widetilde{g}}]$ where $\widetilde{g}$ is a preimage of $g$ in $\Aut(\calH)$.
	We remark that $\overline{\alpha^{\widetilde{g}}}=\overline{\alpha_i}^g$ for $\alpha\in \calH$.
	We may assume that $\calB(v^g,w^g)=\lambda\calB(v,w)^\sigma$ for $v,w\in V$ where some $\lambda\in\bbF_q$ and $\sigma\in\Aut(\bbF_q)$.
	Hence, the image of $w^g=\prod_{i=1}^k[\alpha_i^{\widetilde{g}},\beta_i^{\widetilde{g}}]$ in $N$ is equal to
	\[\prod_{i=1}^kZ\left(\calB\left(\overline{\alpha_i^{\widetilde{g}}},\overline{\beta_i^{\widetilde{g}}}\right)\right)
	=Z\left(\sum_{i=1}^k\calB\left(\overline{\alpha_i}^{g},\overline{\beta_i}^g\right)\right)
	=Z\left(\lambda\left(\sum_{i=1}^k\calB(\overline{\alpha_i},\overline{\beta_i})\right)^\sigma\right)
	=Z(0).\]
	Thus, $w^g\in W$ for any $w\in W$ and $g\in G$.
	
	Hence, $G$ is a subgroup of $\Aut(N)^{V}$.
	By Lemmas~\ref{lem:spext}\,(i) and~\ref{lem:split}, the preimage $A$ of $G$ in $\Aut(N)$ is isomorphic to $p^{2mn^2}{:}\CGammaSp(2m,q)$, and $A^V=G\cong\CGammaSp(2m,q)$.

	It is clear that $V$ is the natural $\bbF_pG$-module.
	We now prove that $\bfZ(N)$ is the natural $\bbF_p\GammaL(1,q)$-module.
	For $g\in G$, we assume that $\calB(v^g,w^g)=\lambda\calB(v,w)^\sigma$ for any $v,w\in V_q$ with $\lambda\in\bbF_q$ and $\sigma\in\Aut(\bbF_q)$.
	Let $g_N$ be the automorphism of $N$ induced by $g$.
	Then $[e_1(a)^g,f_1(1)^g]$ is a preimage of $Z(a)^{g_N}=[E_1(a),F_1(1)]^{g_N}$ in $\calH$.
	Hence, we have
	\[Z(a)^g=Z\left(\calB\left(\overline{e_1(a)}^g,\overline{f_1(1)}^g\right)\right)=Z(\lambda a^\sigma).\]
	This yields that the image of $G$ acting on $\bfZ(N)$ is isomorphic to $\GammaL(1,q)$.
	Thus, we have $A^{\bfZ(N)}=G^{\bfZ(N)}\cong\GammaL(1,q)$.
\end{proof}

\section{Exterior Square of Transitive Linear Groups}\label{sec:mod}

Let $p$ be an odd prime, $V$ be a vector space over $\bbF_p$.
To classify special $3$-orbit $p$-group, Lemma~\ref{lem:spext}\,(iii) and (iv) lead us to focus on the pairs $(G,W)$ such that $W<\Lambda_{\bbF_p}^2(V)$ is invariant under $G\leqs\GL(V)$ and $G$ acts transitively on non-zero vectors of both $V$ and $\Lambda_{\bbF_p}^2(V)/W$.
It was shown in~\cite[Theorem 1.2]{li2024proof} that if $p=2$ and $\dim V=\dim \Lambda_{\bbF_p}^2(V)/W$, then either $G$ is solvable or $G^{(\infty)}\cong \SL_3(2^f)$ for some integer $f$.
This result was applied to prove Gross' conjecture on $2$-automorphic $2$-groups.
For classifying $3$-orbit $p$-group of odd prime $p$, we will prove Theorem~\ref{thm:transmod} in this section, which deals with the case where $p$ is odd and $\dim V$ has not to be equal to $\dim \Lambda_{\bbF_p}^2(V)/W$.

Let $V$ be an irreducible $\F G$-module for some $G\leqs\GL(V)$ over finite field $\F$.
Then $\E=\End_{\bbF_pG}(V)$ is an extension field of $\F$.
So $V$ can be regarded as an $\E G$-module, denoted by $U$ (see~\cite[Section 9]{aschbacher2000Finite}).
Let $\rho: G\rightarrow\GL(U)$ be the homomorphism corresponds to the $\E G$-module $U$.
For $\phi\in\Gal(\E/\F)$, denote by $U^\phi$ the $\E G$-module such that $u^g=(u^{\rho(g)})^\phi$ for $u\in U$ and $g\in G$.
We have the following decomposition given in~\cite[Lemma~2.2]{li2024proof}:
\begin{equation}\label{eq:decomp}
	\Lambda^2_\F(V)=A\oplus B \mbox{ with }B=\bigoplus_{j=1}^{\lfloor \frac{f}{2}\rfloor}B_j\mbox{, where}
\end{equation}
\begin{enumerate}
	\item[(a)] $A\cong \Lambda_\E^2(U)$ as $\F G$-modules;
	\item[(b)] $B_j\cong U\otimes_{\E} U^{\phi^j}$ as $\F G$-modules, when $1\leqs j < f/2$;
	\item[(c)] $2\cdot B_{f/2}\cong U\otimes_{\E} U^{\phi^{f/2}}$ as $\F G$-modules, when $f$ is even.
\end{enumerate}

We remark that a linear group $G\leqs \GL(V)$ acts transitively on non-zero vectors is called a \textit{transitive linear group}.
The classification of transitive linear groups is a part of well-known classification of $2$-transitive groups, see~\cite[Table~7.3]{cameron1999Permutation} for instance.
\begin{theorem}\label{thm:translinear}
	Let $G\leqs\GL(n,p)$ be a transitive linear group with odd prime $p$.
	Then one of the following statements holds.
	\begin{enumerate}
		\item $G$ is solvable.
		\item $G^{(\infty)}$ is isomorphic to $\SL(m,p^f)$ for $m\geqs 2$ and $(m,p^f)\neq (2,3)$; or isomorphic to $\Sp(2m,p^f)$ for $m\geqs 2$.
		\item $G^{(\infty)}\cong \SL(2,5)$ with $(n,p)\in \{(2,11),(2,19),(2,29),(2,59)\}$.
		\item $G^{(\infty)}$ is isomorphic to $2^{1+4}.\A_5$ for $(n,p)=(4,3)$; or isomorphic to $\SL(2,13)$ for $(n,p)=(6,3)$.
	\end{enumerate}
\end{theorem}

First, we consider the case where $\SL(m,q)\lhd G$.
This case is similar to~\cite[Theorem 1.2]{li2024proof}, and we give a simple proof below for completeness.
\begin{lemma}\label{lem:sl1}
	Let $V=\bbF_p^n$, and let $G\leqs\GL(V)$ such that $G^{(\infty)}\cong \SL(m,q)$ with $q=p^f$ odd, $n=mf$ and $m\geqs 3$.
	Then $G$ acts transitively on non-zero vectors of $M=\Lambda^2_{\bbF_p}(V)/W$ if and only if $m=3$ and $M$ is isomorphic to the dual of $V$ as $\bbF_p G^{(\infty)}$-modules. 
	In particular, $W$ is uniquely determined by $G^{(\infty)}$.
\end{lemma}
\begin{proof}
	Note that the induced linear group of $G^{(\infty)}$ on $M$ is either trivial or has a composition factor $\PSL(m,q)$.
	By~\eqref{eq:decomp}, we can decompose $\Lambda^2_{\bbF_p}(V)$ into a direct sum of $\bbF_pG^{(\infty)}$-submodules:
	\[\Lambda^2_{\bbF_p}(V)=A\oplus B=A\oplus\left(\bigoplus_{j=1}^{\lfloor \frac{f}{2}\rfloor}B_j\right).\]
	Let $U$ be the natural $\bbF_q G^{(\infty)}$-module.
	Then $A$, $\bbF_p^{(\infty)}$-isomorphic to $\Lambda^2_{\bbF_q}(U)$, is an irreducible $\bbF_pG^{(\infty)}$-module, see~\cite[Theorem 2.2]{liebeck1987affine}.
	By Steinberg's Twisted Tensor Product Theorem~\cite{steinberg1963Representations}, each of $B_j$ is an irreducible $\bbF_pG^{(\infty)}$-module.
	Moreover, we have that $M$ cannot be a trivial $\bbF_pG^{(\infty)}$-module.
	Thus, the induced linear group of $G^{(\infty)}$ has a composition factor $\PSL(m,q)$.
	By Theorem~\ref{thm:translinear}, we have that $G^{(\infty)}$ acts faithfully and transitively on non-zero vectors of $M$ with $\dim M=\dim V=mf$.
	Then one of $A$ and $B_j$ for $j=1,...,\lfloor\frac{f}{2}\rfloor$ has dimension $n$.
	Remark that 
	\[\begin{aligned}
		&\dim A=f(\dim_{\bbF_q}\Lambda^2_{\bbF_q}(U))=\frac{m(m-1)}{2}f,\\
		&\dim B_j=f(\dim_{\bbF_q} U)^2=m^2f\mbox{ for $1\leqs j<\frac{f}{2}$ and }\\
		&\dim B_{\frac{f}{2}}=\frac{f}{2}(\dim_{\bbF_q} U)^2=\frac{m^2f}{2}\mbox{ when $f$ is even}.
	\end{aligned}\]
	Hence, we have that $M\cong A$ as $\bbF_pG^{(\infty)}$-modules and $m=3$.
	By~\cite[Theorem 2.2]{liebeck1987affine}, $M$ is isomorphic to the dual of $V$ as $\bbF_pG^{(\infty)}$-modules.

	Conversely, if $M$ is isomorphic to the dual of $U$ as $\bbF_pG^{(\infty)}$-modules when $m=3$.
	Then $\dim M=\dim V$, and hence the above arguments imply that $B=W$.
	It suffices to show that $B$ is an $\bbF_pG$-submodule.
	Since $G^{(\infty)}\lhd G$, $B^g$ is an $\bbF_pG_{(\infty)}$-submodule for $g\in G$.
	Suppose that $B^g\neq B$.
	Then $B^g/(B^g\cap B)\cong (B+B^g)/B\cong A$ as $\bbF_pG^{(\infty)}$-modules. 
	As $B^g$ is a completely reducible $\bbF_pG^{(\infty)}$-module, it contains an irreducible $\bbF_pG^{(\infty)}$-submodule $B_0$ such that $\dim B_0=\dim A=mf$.
	This is impossible since $B_0^{g^{-1}}$ is an irreducible $\bbF_pG$-submodule of $B$.
	Thus, we have that $B$ is an $\bbF_pG^{(\infty)}$-submodule for $g\in G$.

	The ``in particular'' part follows since $W$ is always equal to $B$, which is determined by $G^{(\infty)}$.
\end{proof}

Now we deal with the case where $\Sp(2m,q)\lhd G$, which is the main difference between Theorem~\ref{thm:transmod} and~\cite[Theorem 1.2]{li2024proof}.
\begin{lemma}\label{lem:sp1}
	Let $V=\bbF_p^n$, and let $G\leqs\GL(V)$ such that $G^{(\infty)}\cong \Sp(2m,q)$ with $q=p^f$ odd, $n=2mf$ and $m\geqs 1$.
	\begin{enumerate}
		\item There exists a unique $\bbF_pG$-submodule $S$ of $\Lambda^2_{\bbF_p}(V)$ with codimension $f$ such that $G^{(\infty)}$ acts trivially on $\Lambda^2_{\bbF_p}(V)/S$.
		\item If $G$ acts transitively on non-zero vectors of $\Lambda^2_{\bbF_p}(V)/W$, then $W\supseteq S$.
		\item If $\dim \Lambda^2_{\bbF_p}(V)/W\leqs f$, then $W\supseteq S$.
	\end{enumerate}
\end{lemma}
\begin{proof}
	Let $U$ be the natural $\bbF_qG^{(\infty)}$-module.
	We decompose $\Lambda^2_{\bbF_p}(V)$ into a direct sum of $\bbF_pG^{(\infty)}$-submodules as in~\eqref{eq:decomp}:
	\[\Lambda^2_{\bbF_p}(V)=A\oplus B=A\oplus\left(\bigoplus_{j=1}^{\lfloor \frac{f}{2}\rfloor}B_j\right).\] 
	By Steinberg's Twisted Tensor Product Theorem, we have that each of $B_j$ for $j=1,...,\lfloor\frac{f}{2}\rfloor$ is a faithful and irreducible $\bbF_pG^{(\infty)}$-module.
	Remark that $A\cong \Lambda_{\bbF_q}^2(U)$ has a $\bbF_qG^{(\infty)}$-submodule $T$ of codimension $1$ over $\bbF_q$ (see~\cite{carter2005Lie,korhonen2019Jordan,mcninch1998Dimensional}) such that
	\begin{enumerate}
		\item[(a)] if $p\nmid m$ then $T$ is an irreducible $\bbF_pG^{(\infty)}$-submodule; 
		\item[(b)] if $p\mid m$ then $T$ is indecomposable such that $T/\mathrm{rad}(T)$ is an irreducible $\bbF_pG^{(\infty)}$-module, where $\dim_{\bbF_q}\mathrm{rad}(T)=1$.
	\end{enumerate}
	Let $R$ be the radical of $\Lambda^2_{\bbF_p}(V)$ as $\bbF_pG$-module.
	Then either $\Lambda^2_{\bbF_p}(V)$ is completely reducible when $p\nmid m$, or $R=T$ when $p\mid m$.

	\textbf{Part~(i).} 
	It is clear that $S\supseteq B$, otherwise $\Lambda_{\bbF_p}^2(V)/S$ is $\bbF_pG^{(\infty)}$-isomorphic to $B_j$ for some $j$.
	Thus, $\Lambda^2_{\bbF_p}(V)/S$ is naturally $\bbF_pG^{(\infty)}$-isomorphic to some quotient of $A$.
	Then we obtain that $S\supseteq T+B$.
	Note that $\dim A/T=f$, and thus $S=T+B$.
	It suffices to show that $S$ is an $\bbF_pG$-submodule.
	Suppose that $S^g\neq S$ for some $g\in G$.
	Then $S^g/(S\cap S^g)\cong (S+S^g)/S$ is a trivial $\bbF_pG^{(\infty)}$-module.
	This yields that $G^{(\infty)}$ acts trivially on $S/(S^{g^{-1}}\cap S)$.
	Remark that $\mathrm{rad}(S)=\mathrm{rad}(T)\oplus B$ contains no trivial $\bbF_pG^{(\infty)}$-submodules, a contradiction.
	Hence, $S$ is an $\bbF_pG$-module.

	\textbf{Part~(ii).} 
	Suppose that $M$ is not a trivial $\bbF_pG^{(\infty)}$-module.
	Let $M=\Lambda^2_{\bbF_p}(V)/W$.
	By Theorem~\ref{thm:translinear}, we have that $G^{(\infty)}$ acts faithfully and transitively on non-zero vectors of $M$ with $\dim M=\dim V$.
	Since $M$ is $\bbF_pG^{(\infty)}$-isomorphic to a submodule of $\Lambda^2_{\bbF_p}(V)/R$, it yields that $M$ is $\bbF_pG^{(\infty)}$-isomorphic to either $B_j$ for some $j$ or $T$ when $p\nmid m$.
	Note that 
	\[\begin{aligned}
		&\dim T=f\left(\frac{2m(2m-1)}{2}-1\right),\\
		&\dim B_j=f(\dim_{\bbF_q} U)^2=4m^2f\mbox{ for $1\leqs j< \frac{f}{2}$ and }\\
		&\dim B_{\frac{f}{2}}=\frac{f}{2}(\dim_{\bbF_q} U)^2=2m^2f\mbox{ when $f$ is even}.
	\end{aligned}\]
	By comparing dimensions, we have that $M$ is $\bbF_pG^{(\infty)}$-isomorphic to $B_{ \frac{f}{2}}$ with $m=1$ and $f$ is even.
	Recall that $2\cdot B_{ \frac{f}{2}}\cong U\otimes U^{\phi^{\frac{f}{2}}}$ and $G^{(\infty)}\cong \SL(2,q)$.
	Then the induced linear group of $G^{(\infty)}$ acting on $2\cdot B_{\frac{f}{2}}$ is isomorphic to $\Omega^{-}_4(p^{f/2})$ (see.~\cite[page 45]{kleidman1990subgroup}), which is not transitive on non-zero vectors.
	Thus, we have that $M$ is a trivial $\bbF_pG^{(\infty)}$-module.

	Note that $B$ is a completely reducible $\bbF_pG^{(\infty)}$ with no trivial $\bbF_p^{(\infty)}$-submodules.
	This implies that $W\supset B$.
	Thus, we have that $S=T\oplus B\subseteq W$ when $p\mid m$ since $T=R$ is the radical of $\Lambda_{\bbF_p}^2(V)$.
	When $p\nmid m$, we have that $T$ is a non-trivial irreducible $\bbF_p^{(\infty)}$-submodule of $A$.
	Thus, $S=T\oplus B\subseteq W$ as desired.

	\textbf{Part~(iii).}
	Note that each non-trivial $\bbF_pG^{(\infty)}$-submodule of $\Lambda^2_{\bbF_p}(V)/R$ has dimension larger than $f$.
	Then $\Lambda^2_{\bbF_p}(V)/W$ is trivial.
	Therefore, $S\subseteq W$ since $S/(S\cap R)$ has no trivial submodules.
\end{proof}

Now we are ready to give a proof for Theorem~\ref{thm:transmod}.
\begin{proof}[Proof of Theorem~\ref{thm:transmod}]
	Note that $G$ is isomorphic one of the groups in Theorem~\ref{thm:translinear}.
	Using \Magma~\cite{magma}, we obtain that if $G$ is one of the groups listed in Theorem~\ref{thm:translinear}\,(iii) and (iv), then $\dim \Lambda^2_{\bbF_p}(V)/W=1$, as in case~(i).
	If $G^{(\infty)}\cong \SL(m,q)$ with  $m\geqs 3$, $q=p^f$ and $n=mf$, then case~(ii) holds by Lemma~\ref{lem:sl1}.
	Assume that $G_0\cong \Sp(2m,p^f)$ with $(m,p^f)\neq(1,3)$ and $n=2mf$.
	Lemma~\ref{lem:sp1} implies case~(iii).
\end{proof}

\section{$3$-orbit $p$-groups for odd prime $p$}\label{sec:oddp}

Let $N$ be a finite non-abelian $3$-orbit $p$-group with $p$ an odd prime.
Assume that $N/\Phi(N)\cong\bbZ_p^n$.
Then $N$ is a quotient group $\calH/W$ where $\calH=\calH_{n,p}$ and $W<\bfZ(\calH)$.
Identify $N/\bfZ(N)$ and $\bfZ(N)$ with $\calH/\bfZ(\calH)$ and $\bfZ(\calH)/W$, respectively.
Let $V=N/\bfZ(N)$ and $M=\bfZ(N)$.
Recall Lemmas~\ref{lem:spext}\,(iii) and~\ref{lem:oddp} that
\begin{enumerate}
	\item $W$ is an $\bbF_pG$-submodule of $\Lambda^2_{\bbF_p}(V)$;
	\item $G=\Aut(N)^V$ is the maximal subgroup of $\Aut(\calH)^V=\GL(V)$ in the sense that $W^G=W$; and
	\item $G$ acts transitively on non-zero vectors of both $V$ and $M\cong \Lambda^2_{\bbF_p}(V)/W$.
\end{enumerate}
We may assume that $G$ is non-solvable due to Theorem~\ref{thm:revdorn}.
If $\dim M=1$, then $N$ is an extraspecial $p$-group, and hence $N\cong p^{1+2m}_{+}$ for some $m\geqs 1$.
Thus, we mainly focus on cases that $G^{(\infty)}$ is isomorphic to $\SL(3,p^f)$ or $\Sp(2m,p^f)$ by Theorem~\ref{thm:transmod}.

\subsection{$\SL(3,q)\lhd G$}\

Assume that $G_0=G^{(\infty)}\cong \SL(3,q)$ with $q=p^f$.
We will show that such $3$-orbit group is uniquely determined by its order, and it is isomorphic to $A_p(n,\theta)$ with $|\theta|=3$ given in Example~\ref{exam:apn}, as in the 6-th row of Table~\ref{tab:3aogp}.

First, we show the existence and the uniqueness of non-abelian $3$-orbit $p$-group $N$ with $|N|=|\bfZ(N)|^2$.
\begin{lemma}\label{lem:slunique}
	Let $p$ be an odd prime.
	Then the following statements hold.
	\begin{enumerate}
		\item If $N$ is a $3$-orbit $p$-group such that $|N|=|\bfZ(N)|^2$, then $|\bfZ(N)|=q^3$ for some $q=p^f$ and $\SL(3,q)\lhd \Aut(N)^{N/\bfZ(N)}$.
		\item When $|N|=q^6$ for some $q=p^f$, there exists a unique $3$-orbit $p$-group such that $|\bfZ(N)|=q^3$.
		\item Let $N$ be the $3$-orbit $p$-group with $|N|=|\bfZ(N)|^2=q^6$ and $q=p^f$.
		Then $\Aut(N)\cong p^{9f^2}{:}\GammaL(3,q)$.
	\end{enumerate}
\end{lemma}
\begin{proof}
	Assume that $N$ is a $3$-orbit $p$-group such that $|N|=|\bfZ(N)|^2$.
	Then $\Aut(N)$ is non-solvable by Theorem~\ref{thm:revdorn}.
	Let $V=N/\bfZ(N)$, $M=\bfZ(N)$ and $G=\Aut(N)^V$.
	Lemma~\ref{lem:spext} yields that $G$ acts transitively on non-zero vectors of $V$ and $M$ and $M\cong \Lambda^2_{\bbF_p}(V)/W$ for some $\bbF_pG$-submodule $W$ of $\Lambda^2_{\bbF_p}(V)$.
	Hence, we have that $G^{(\infty)}\cong \SL(3,q)$ with $|M|=|V|=q^3$ by Theorem~\ref{thm:transmod}, as in part~(i).
	
	Let $\calH=\calH_{n,p}$ with $n=3f$, and let $V=\calH/\bfZ(\calH)$.
	Then $X=\Aut(\calH)^V=\GL(V)\cong\GL(n,p)$ and $\bfZ(\calH)\cong\Lambda^2_{\bbF_p}(V)$ as $\bbF_pX$-modules by Lemma~\ref{lem:oddp}\,(i).
	Let $G\leqs X$ such that $G\cong\GammaL(3,q)$.
	Lemma~\ref{lem:sl1} shows that there exists an $\bbF_pG$-submodule $W$ of $\Lambda^2_{\bbF_p}(V)$ such that $\Lambda^2_{\bbF_p}(V)/W$ is $\bbF_pG^{(\infty)}$-isomorphic to the dual of $V$.
	We identify $W$ as a subgroup of $\bfZ(\calH)$.
	Then $G\leqs \Aut(N)^V$ by Lemma~\ref{lem:oddp}\,(ii).
	It follows that $G$ acts transitively on non-identity elements of both $N/\bfZ(N)$ and $\bfZ(N)$.
	Hence, $N$ is a $3$-orbit group by Lemma~\ref{lem:spext}\,(iii).
	This proves the existence in part~(ii).
	Now we prove the uniqueness in part~(ii) by assuming that $N_1=\calH/W_1$ is another desired $3$-orbit $p$-group.
	Then part~(i) implies that $\SL(3,q)\lhd \Aut(N_1)^V$.
	Denote $\Aut(N_1)$ by $G_1$.
	Note that there exists $x\in X$ such that $(G_1)^{(\infty)}$ is conjugate to $G^{(\infty)}$ by $x$, and hence $G_1^x\leqs G$.
	Note that $W_1^x$ and $W$ are $\bbF_pG_1^x$-submodules such that $G_1^x$ acts transitively on the non-zero vectors of the corresponding factor modules.
	Lemma~\ref{lem:sl1} yields that $W_1^x=W$.
	By \ref{lem:oddpiso}, we have that 
	\[N_1=\calH/W_1\cong \calH/W_1^x=\calH/W\cong N.\]
	
	Let $N=\calH/W$ as defined above.
	Then $\GammaL(3,q)\cong G\leqs\Aut(N)^V$ as given in the proof of the existence in part~(ii).
	By part~(i), we have that $\SL(3,q)\lhd\Aut(N)^V$.
	It follows that $G=\Aut(N)^V$, and hence $\Aut(N)\cong p^{n^2}{:}\GammaL(3,q)$ by Lemmas~\ref{lem:spext}\,(i)~and~\ref{lem:split}, as in part~(iii).
\end{proof}

Recall that $|A_p(n,\theta)|=|\bfZ(A_p(n,\theta))|^2=p^{2n}$ given in Example~\ref{exam:apn}.
We aim to prove $A_p(n,\theta)$ is a $3$-orbit group when $3\mid n$ and $|\theta|=3$.
Observe that $\GammaL(3,p^{n/3})$ contains a Singer cycle, the cyclic subgroup of order $p^n-1$, and $\Aut(A_p(n,\theta))$ also has an automorphism of order $p^n-1$ given in Proposition~\ref{prop:antheta}\,(ii).
We need the following lemma before proving the isomorphism.
\begin{lemma}\label{lem:cyclicexterior}
	Suppose that $V=\bbF_p^n$.
	Let $G\leqs\GL(V)\cong\GL(n,p)$ be the Singer cycle, the cyclic subgroup of $\GL(V)$ of order $p^n-1$.
	Then $\Lambda^2_{\bbF_p}(V)$ is a direct sum of pairwise non-isomorphic irreducible $\bbF_p G$-submodules.
\end{lemma}
\begin{proof}
	Let $\E=\bbF_q$ with $q=p^n-1$, and let $U$ be the natural $\bbF_q G$-module.
	Clear that $\End_{\bbF_pG}(V)\cong\E$.
	We may assume that $U\cong V$ as $\bbF_pG$-modules. 
	Notice that $\Lambda_{\bbF_qG}^2(U)=1$.
	Let $\phi$ be the Frobenius automorphism $x\mapsto x^p$ in $\E$.
	With the decomposition given in~\eqref{eq:decomp}, we have that 
	\[\Lambda^2_{\bbF_p}(V)\cong \bigoplus_{i=1}^{\lfloor \frac{n}{2}\rfloor} B_j,\]
	where, as $\bbF_pG$-modules, $B_j\cong U\otimes_{\E} U^{\phi^j}$ for $1\leqs j<\frac{n}{2}$; and $2\cdot B_{\frac{n}{2}}\cong U\otimes_{\E} U^{\phi^{n/2}}$.
	It suffices to show that each $B_j$ is an irreducible $\bbF_p G$-module and $B_i\ncong B_j$ if $i\neq j$.

	First, we show that each $B_j$ is irreducible.
	Note that the kernel of $G$ on $B_j$ has order $\gcd(p^j+1,p^n-1)$.
	Then, when $n$ is even, the permutation image of $G$ on $B_{\frac{n}{2}}$ is isomorphic to $\bbZ_{p^{\frac{n}{2}}-1}$, which is irreducible on $B_{\frac{n}{2}}$.
	Assume that $j<\frac{n}{2}$ with $n\geqs 3$.
	If $(n,p)=(6,2)$, it is not hard to check that $B_1$ and $B_2$ are non-isomorphic irreducible $\F G$-modules.
	If $(n,p)\neq (6,2)$ and $n\geqs 3$, there exists a prime divisor $r$ of $p^n-1$ such that $r\mid p^n-1$ and $r\nmid p^k-1$ for any $k<n$ by Zsigmondy's Theorem~\cite{zsigmondy1892Zur}.
	Then the order of the induced permutation group of $G$ on $B_j$ is divisible by $r$, and hence $B_j$ is an irreducible $\bbF_pG$-module as $\dim B_j=n$.

	It is clear that $B_{\frac{n}{2}}\ncong B_j$ for $1\leqs j<\frac{n}{2}$ as they have different dimensions.
	We only need to show that $B_i\ncong B_j$ for $1\leqs i\neq j<n/2$.
	Note that $B_i\cong B_j$ if and only if $U\otimes_\E U^{\phi^i}\cong (U\otimes_\E U^{\phi^j})^{\phi^k}$ for some $k$.
	Let $g$ be a generator of $G$, and let $\lambda$ be the character of $g$ on $U$.
	Suppose that $U\otimes_\E U^{\phi^i}\cong (U\otimes_\E U^{\phi^j})^{\phi^k}$.
	Then 
	\[\lambda^{1+p^i}=\lambda^{p^k(1+p^i)}.\]
	Note that $\lambda$ is a generator of $\E^\times$.
	Hence, we have that
	\[1+p^i\equiv p^k+p^{j+k}\pmod{p^n-1}.\]
	Let $\ell$ be the residue of $j+k$ modulo $n$.
	Then
	\[1+p^i\equiv p^k+p^{\ell}\pmod{p^n-1}.\]
	Note that $0\leqs i,k,\ell <n$, we have
	\[-(p^n-1)<-(p^{i}-1)\leqs p^{k}+p^{\ell}-p^{i}-1<p^n-1.\]
	Hence, the characters of $g$ on $U\otimes_\E U^{\phi^i}$ and $(U\otimes_\E U^{\phi^i})^{\phi^k}$ cannot be equal.
	Therefore, $\Lambda^2_{\bbF_p}(V)$ is a direct sum of non-isomorphic irreducible $\bbF_pG$-modules.
\end{proof}

Now we are ready to get the 6-th row of Table~\ref{tab:3aogp}.
\begin{lemma}\label{lem:autoapsl}
	Let $N=A_p(n,\theta)$ be defined in Example~\ref{exam:apn} such that $3\mid n$, $|\theta|=3$ and $p$ is an odd prime.
	Then $N$ is a $3$-orbit group and $\Aut(N)\cong p^{n^2}{:}\GammaL(3,p^{n/3})$.
\end{lemma}
\begin{proof}
	Note that $N$ is a special $p$-group of exponent $p$.
	We assume that $N=\calH/W$ with $\calH=\calH_{n,p}$ and $W<\bfZ(\calH)$.
	Let $V=\calH/\bfZ(\calH)=N/\bfZ(N)$, and let $X=\Aut(\calH)^V=\GL(V)\cong\GL(n,p)$.
	Recall that the automorphism $\xi\in\Aut(N)$, given in Proposition~\ref{prop:antheta}\,(ii), has order $p^n-1$ and $\langle\xi\rangle$ acts faithfully on $V$.
	Hence, $\langle \xi\rangle$ induces a subgroup $H\leqs X$ of order $p^n-1$ and $W^H=W$ by Lemma~\ref{lem:oddp}\,(ii).
	Then there exists $G\leqs X$ such that $H<G$ and $G\cong \GL(3,q)$ with $q=p^{n/3}$.
	By Lemma~\ref{lem:sl1}, there exists a unique $\bbF_pG$-submodule $W_1<\Lambda^2_{\bbF_p}(V)$ such that $G$ acts transitively on non-zero vectors of $\Lambda^2_{\bbF_p}(V)/W_1$.
	We identify $W_1$ as a subgroup of $\bfZ(\calH)$.
	Then $N_1=\calH/W_1$ is a $3$-orbit group as given in the proof of Lemma~\ref{lem:slunique}\,(i).
	It suffices to show that $W_1=W$, and then $N=N_1$ is a $3$-orbit group with $\Aut(N)\cong p^n{:}\GammaL(3,p^{n/3})$ as given in Lemma~\ref{lem:slunique}\,(iii).

	Remark that each of $\Lambda_{\bbF_p}(V)$, $W$ and $W_1$ is an $\bbF_pH$-module and is a direct sum of non-isomorphic irreducible $\bbF_pH$-submodules by Lemma~\ref{lem:cyclicexterior}.
	Thus, we have that $W_1=W$ if and only if $\Lambda_{\bbF_p}(V)/W\cong \Lambda_{\bbF_p}(V)/W_1$ as $\bbF_pH$-modules.
	It is clear that both of these two $\bbF_pH$-modules are irreducible.
	Let $\chi$ and $\chi_1$ be the characters of $\Lambda_{\bbF_p}(V)/W$ and $\Lambda_{\bbF_p}(V)/W_1$.
	We only need to verify that $\chi(\xi)=\chi_1(\xi)$.
	Let $U$ be the natural $\bbF_qH$-module with $q=p^n$.
	Assume that $U\cong V$ as $\bbF_pH$-modules.
	Let $\lambda$ be the $\bbF_q$-character of $\xi$ on $U$.
	Then $\lambda$ generates $\bbF_q^\times$ and the character of $\xi$ acting on $V$ is $\mathrm{Tr}_{\bbF_{q}/\bbF_p}(\lambda)$.
	By definition of $\xi$, we have that
	\[\chi(\xi)=\mathrm{Tr}_{\bbF_{q}/\bbF_p}(\lambda\cdot \lambda^\theta)=\sum_{i=1}^n\lambda^{p^i(1+p^{\frac{n}{3}})}.\]
	Recall that $G\cong \GL(3,p^{n/3})$.
	Let $\mathcal{U}$ be the natural $\bbF_{p^{\frac{n}{3}}}G$-module.
	Let $\chi_{\mathcal{U}}$ be the $\bbF_{\frac{n}{3}}$-character of $G$ on $\mathcal{U}$.
	Then 
	\[\chi_{\mathcal{U}}(\xi)=\mathrm{Tr}_{\bbF_q/\bbF_{p^{n/3}}}(\lambda)=\lambda+\lambda^{p^{\frac{n}{3}}}+\lambda^{p^{\frac{2n}{3}}}.\]
	In the proof of Lemma~\ref{lem:sl1}, we obtain that $\Lambda_{\bbF_p}(V)/W_1$ is $\bbF_{p}G$-isomorphic to $\Lambda_{\bbF_{\frac{n}{3}}}^2(\mathcal{U})$.
	Thus, we have that 
	\[\begin{aligned}
		\chi_1(\xi)&=\mathrm{Tr}_{\bbF_{p^{n/3}}/\bbF_p}\left( \frac{\chi_{\mathcal{U}}(\xi)^2-\chi_{\mathcal{U}}(\xi^2)}{2}\right)=\sum_{i=1}^{\frac{n}{3}}\left( \lambda^{1+p^{\frac{n}{3}}}+\lambda^{1+p^{\frac{2n}{3}}}+\lambda^{p^{\frac{n}{3}}+p^{\frac{2n}{3}}}\right)^{p^i}\\
		&=\sum_{i=1}^{\frac{n}{3}}\left( \sum_{k=1}^3\left(\lambda^{p^{\frac{kn}{3}}(1+p^{\frac{n}{3}})}\right)\right)^{p^i}=\sum_{i=1}^{n}\lambda^{p^i(1+p^{\frac{n}{3}})}=\chi(\xi).
	\end{aligned}\]
	This yields that $W=W_1$, and therefore $N=N_1$ is a $3$-orbit group with $\Aut(N)^{N/\bfZ(N)}\cong\GammaL(3,p^{\frac{n}{3}})$.
\end{proof}

\subsection{$\Sp(2m,p^n)\lhd G$}\ 

Now we treat the case where $\Sp(2m,q)\lhd G$.
The following lemma shows that these groups are factor groups of $q^{1+2m}_{+}$.
\begin{lemma}\label{lem:nmpncgamma}
	Let $q=p^n$ with $p$ an odd prime, and let $N$ be a special $p$-group of exponent $p$ with $N/\bfZ(N)\cong\bbZ_p^{2mn}$ for some integer $m$.
	The following statements hold.
	\begin{enumerate}
		\item If $|\bfZ(N)|= q$ and $\Aut(N)^{N/\bfZ(N)}$ has a subgroup isomorphic to $\Sp(2m,q)$, then $N\cong q^{1+2m}_{+}$.
		\item If $|\bfZ(N)|\leqs q$ and $\Aut(N)^{N/\bfZ(N)}$ has a subgroup isomorphic to $\Sp(2m,q)$, then $N$ is a quotient of $q^{1+2m}_{+}$.
		\item $\Aut(q^{1+2m}_+)\cong p^{2mn^2}{:}\CGammaSp(2m,q)$.
	\end{enumerate}
\end{lemma}
\begin{proof}
	Let $\calH=\calH_{2mn}$, and let $V=\calH/\bfZ(\calH)$.
	Then $X=\Aut(\calH)^{V}=\GL(V)\cong\GL(2mn,p)$ and $\bfZ(\calH)$ can be identified with $\Lambda^2_{\bbF_p}(V)$ as $\bbF_pX$-modules by Lemma~\ref{lem:oddp}\,(i).
	Note that $N=\calH /W$ for some $W<\bfZ(N)$.
	Let $G=\Aut(N)^V$.
	Then $W^g=W$ for any $g\in G$.
	
	Suppose that $G_0\leqs G$ isomorphic to $\Sp(2m,q)$ and $|\bfZ(N)|=|\Lambda^2_{\bbF_p}(V)/W|=q$.
	Note that $q^{1+2m}_+=\calH/W_0$ for some $W\leqs\bfZ(\calH)$ and $X_{W_0}$ contains a subgroup $G_1$ isomorphic to $\Sp(2m,q)$ by Lemma~\ref{lem:oddp}\,(ii)~and~\ref{lem:autnmq}.
	Then there exists $x\in X$ such that $G_1^x=G_0$, and then $q^{1+2m}_{+}\cong\calH/W^x$ with $G_0\leqs X_{W^x}$.
	Thus, we may assume that $W_0^x=W_0$ and $G_1=G_0$.
	It suffices to show that $W_0=W$.
	By Lemma~\ref{lem:sp1}\,(iii), $S\subseteq W$ and $S\subseteq W_0$, where $S$ is unique determined by $G_0=G^{(\infty)}$.
	Note that $\dim S=\dim \Lambda^2{\bbF_q}(V)-n$, it follows that $S=W=W_0$ as $|\bfZ(\calH)/S|=|\bfZ(N)|=|\bfZ(q^{1+2m}_+)|=p^n$.
	Thus, $N\cong q^{1+2m}_+$, as in part~(i).

	Part~(ii) can be proved with similar arguments as above by Lemma~\ref{lem:sp1}\,(i).

	Let $N=q^{1+2m}_+$.
	If $q=p$, then $\Aut(N)\cong p^{2m}{:}\mathrm{CSp}(2m,q)$ given by Winter~\cite{winter1972automorphism}.
	Suppose that $n\geqs 2$.
	Then $G=\Aut(N)^V$ contains a subgroup isomorphic to $\CGammaSp(2m,q)$ by Lemma~\ref{lem:autnmq}.
	Hence, $G$ is non-solvable.
	Remark that $G$ acts transitively on non-zero vectors of both $V$ and $\Lambda^2_{\bbF_p}(V)/W$ and $|V|\neq |\Lambda^2_{\bbF_p}(V)/W|$.
	Then $G$ has a normal subgroup isomorphic to $\Sp(2m',p^{n'})$ with $2m'n'=2mn$ by Theorem~\ref{thm:transmod}.
	As $\CGammaSp(2m,q)\lesssim G$, we have that $n'\leqs n$.
	Suppose that $n'<n$.
	Then $|\Lambda^2_{\bbF_p}(V)/W|<q$ by Theorem~\ref{thm:transmod}, a contradiction.
	Thus, $G$ has a normal subgroup isomorphic to $\Sp(2m,q)$.
	Therefore, $G\cong\CGammaSp(2m,q)$ and we obtain part~(iii) by Lemmas~\ref{lem:spext}\,(i)~and~\ref{lem:split}.
\end{proof}

We remark that the $1$-subspace stabilizer $G_\alpha$ of $G=\Sp(2m+2,q)$ has form 
\[q^{1+2m}{:}(\Sp(2m,q)\times\GL(1,q)).\]
Let $N=\textbf{O}_p(G_\alpha)\cong q^{1+2m}$.
Then $N$ is a special $p$-group with $|\bfZ(N)|=q$ and $\Sp(2m,q)\lesssim \Aut(N)^{N/\bfZ(N)}$ (see~\cite[Section 3.5.4]{wilson2009finite}).
These facts together with Lemma~\ref{lem:nmpncgamma}\,(i) yield that $N\cong q^{1+2m}_{+}$, which is Theorem~\ref{thm:isos}\,(iii).
\begin{corollary}\label{coro:psp}
	Let $G=\PSp(2m+2,q)$ where $q=p^n$ with odd prime $p$. 
	Then $G$ naturally acts on $1$-spaces of $\bbF_{q}^{2m+2}$ and $\textbf{O}_p(G_\alpha)\cong q^{1+2m}_{+}$ for any stabilizer $G_\alpha$.
\end{corollary}

\section{Subfield hyperplanes and quotients of $q^{1+2m}_+$}\label{sec:quoextra}

Let $q_0=p^d$ and $q=p^n$ with $d\mid n$.
Note that $\Sp(2m,q)\leqs\Sp(2mn/d,q_0)$ for any positive integer $m$.
Lemma~\ref{lem:nmpncgamma}\,(ii) implies that there exists $U<\bfZ(q^{1+2m}_+)$ such that
\[(q_0)^{1+2mn/d}_{+}\cong q^{1+2m}_{+}/U.\]
To describe the subgroup $U$, we introduce the following concept of \textit{subfield hyperplanes}.
\begin{definition}\label{def:subfield}
	Let $V=\bbF_q^+$ with $q=p^n$ for some prime $p$, and let $x\in H$ be the linear transformation given by multiplying a generator $\lambda$ of $\bbF_q^\times$.
	For any divisor $d$ of $n$, we say an $\bbF_p$-subspace $U$ of $V$ is a \textit{subfield hyperplane} with respect to subfield $\bbF_{p^d}$ if $\dim U=n-d$ and $U$ is invariant under $\langle x^\ell\rangle$ with $\ell=\frac{q-1}{p^d-1}$.  
\end{definition}

Note that $V=\bbF_{p^n}^+$ is naturally an $\bbF_{p^d}$-space for any divisor $d$ of $n$.
Hence, an $\bbF_p$-subspace $U$ of $V$ is a subfield hyperplane with respect to $\bbF_{p^d}$ if and only if $U$ is naturally a hyperplane of $V$ when viewing both $U$ and $V$ as $\bbF_{p^d}$-spaces. 
We provide the following lemma before studying $q^{1+2m}_{+}/U$.
\begin{lemma}\label{lem:subequal}
	Suppose that $H=\GammaL(1,p^n)$ acting naturally on $V=\bbF_{p^n}^+$.
	For any divisor $d$ of $n$, the set of subfield hyperplane with respect to subfield $\bbF_{p^d}$ of $V$ forms an orbit of $H$.

	In particular, $(H_U)^{V/U}\cong \GammaL(1,p^d)$ for any subfield hyperplane $U$ with respect to $\bbF_{p^d}$.
\end{lemma}
\begin{proof}
	Let $\Sigma$ be the set of subfield hyperplanes with respect to $\bbF_{p^d}$.
	Set $x\in H$ the operator admitted by multiplying a generator $\lambda\in\bbF_{p^n}^\times$, and set $y$ the operator admitted by the Frobenius automorphism in $\Aut(\bbF_{p^n})$.

	By definition, $U\in\Sigma$ if and only if $U$ is an $\bbF_p\langle x^\ell\rangle$-submodule of $V$ where $\ell=\frac{p^n-1}{p^d-1}$.
	Since $\langle x^{\ell}\rangle\lhd H$, $U^h$ is also an $\bbF_p\langle x^\ell\rangle$-submodule for any $h\in H$.
	Thus, $\Sigma$ is closed under the action of $H$.

	Note that each $U\in \Sigma$ is a hyperplane of $\bbF_{p^d}^{n/d}$, and hence $|\Sigma|=\frac{p^n-1}{p^d-1}$.
	Since $\langle x\rangle_U=\langle x^\ell\rangle$, it follows that 
	\[|U^{\langle x\rangle}|=\frac{|x|}{|x^{\ell}|}=\frac{p^n-1}{p^d-1}=|\Sigma|.\]
	Thus, $\langle x\rangle$ acts transitively on $\Sigma$, so does $H$.

	Let $U\subseteq V=\bbF_{p^n}^+$ consisting of elements $u$ such that 
	\[\Tr_{\bbF_{p^n}/\bbF_{p^d}}(u)=\sum_{k=1}^{n/d} u^{y^{dk}}=0.\]
	Clearly, $U$ is a hyperplane of $\bbF_{p^d}^{n/d}$ as $\bbF_{p^d}$-space.
	Thus, $U\in\Sigma$ and $U^{x^\ell}=U$.
	In particular, we have 
	\[\Tr_{\bbF_{p^n}/\bbF_{p^d}}(u^y)=\Tr_{\bbF_{p^n}/\bbF_{p^d}}(u)^y=0.\]
	Hence, $U^y=U$ and $U$ is an $\bbF_p\langle x^\ell,y\rangle$-submodule of $V$.
	Since $\langle x^\ell,y\rangle^{V/U}\cong\GammaL(1,p^d)$, we complete the proof.
\end{proof}

Let $N=q^{1+2m}_+$.
Then $\Aut(N)^{\bfZ(N)}\cong\GammaL(1,q)$ by Lemmas~\ref{lem:autnmq} and~\ref{lem:nmpncgamma}(iii).
Hence, $\bfZ(N)$ can be naturally viewed as $\bbF_q^+$ by the action of $\GammaL(1,q)$.
We are now ready to complete the proof of Theorem~\ref{thm:isos}.

\begin{proof}[Proof of Theorem~\ref{thm:isos}]
	Let $N=q^{1+2m}_+$, and let $N_0=(q_{0})^{1+2mn/d}_+$.
	By Lemma~\ref{lem:nmpncgamma}\,(iii), we have $\Aut(N_0)^{N_0/\bfZ(N_0)}\cong\CGammaSp(2mn/d,q_0)$ which has a subgroup isomorphic to $\Sp(2m,q)$.
	Hence, $N_0\cong N/U_0$ for some $U_0<\bfZ(N)$ by Lemma~\ref{lem:nmpncgamma}\,(i).
	Set $H=\Aut(N)^{\bfZ(N)}$.
	Then $H\cong\GammaL(1,q)$ by Lemma~\ref{lem:autnmq}.
	Let $x$ be an element of order $q-1$ in $H$, and let $\delta$ be a preimage of $x$ in $\Aut(N)^{N/\bfZ(N)}$.
	Then $\delta^\ell\in \Aut(N_0)^{N_0/\bfZ(N_0)}\cong \CGammaSp(2mn/d,q_0)$ for $\ell=\frac{q-1}{q_0-1}$ (see~\cite[page 112]{kleidman1990subgroup}).
	Hence, $U_0$ is invariant under $\delta^\ell$, equivalently, $U_0^{x^\ell}=U_0$.
	Thus, $U_0<\bfZ(N)$ is a subfield hyperplane of with respect to $\bbF_{p^d}$ by definition.

	Conversely, assume that $U<\bfZ(N)$ is a subfield hyperplane with respect to $\bbF_{p^d}$.
	Then there exists an $h\in H\cong\GammaL(1,q)$ such that $U^h=U_0$ by Lemma~\ref{lem:subequal}.
	Applying Lemma~\ref{lem:oddpiso} by getting a preimage of $h$ in $\Aut(\calH_{2mn,p})$, we can obtain that 
	\[q^{1+2m}_+/U\cong q^{1+2m}_+/U_0\cong (q_{0})^{1+2mn/d}_+.\]
	Notice that any subspace of codimension $1$ is a subfield hyperplane with respect to $\bbF_{p}$, then the ``in particular'' part immediately holds.
	The proof of part~(i) is completed.

	Additionally, parts~(ii) and~(iii) have been given in Lemma~\ref{lem:susl} and Corollary~\ref{coro:psp}, respectively.
\end{proof}

Let $N=q^{1+2m}_+$ with $q=p^n$, and let $U<\bfZ(N)$.
Suppose that $U$ contains a subfield hyperplane $U_0$ of $\bfZ(N)$ with respect to $\bbF_{q_0}$ with $q_0=p^d$.
By Theorem~\ref{thm:isos}\,(i), we have
\[N/U\cong (N/U_0)/(U_0/U)\cong (q_0)^{1+2mn/d}_+/U_1,\]
where $U_1=U_0/U$.
Thus, we may focus on the case where $U$ contains no subfield hyperplanes, as in the following example.
\begin{example}\label{exam:quot}
	For $q=p^n$ with odd prime $p$, let $N=q^{1+2m}_+/U$ such that
	\begin{enumerate}
		\item $U<\bfZ(q^{1+2m}_+)\cong\bbF_q^+$ contains no subfield hyperplanes with respect to any subfield; 
		\item $\GammaL(1,q)_U$ is transitive on $\bfZ(N)=\bfZ(q^{1+2m}_+)/U$.
	\end{enumerate}
\end{example}

The groups in the above example form the last row in Table~\ref{tab:3aogp}.
\begin{lemma}\label{lem:quotaut}
	Let $N=q^{1+2m}_+/U$ be defined in Example~\ref{exam:quot}.
	Then $N$ is $3$-orbit group and 
	\[\Aut(N)\cong p^{2mn_0n}{:}(\Sp(2m,q){:}\GammaL(1,q)_U),\]
	where $p^{n_0}=|\bfZ(N)|=q/|U|$.
\end{lemma}
\begin{proof}
	Let $V=q^{1+2m}_+/\bfZ(q^{1+2m}_+)$, and identify $N/\bfZ(N)$ with $V$.
	By Lemma~\ref{lem:nmpncgamma}\,(iv), we have $G=\Aut(q^{1+2m}_+)^{V}\cong\CGammaSp(2m,q)$.
	Let $G_0$ be the normal subgroup of $G$ isomorphic to $\Sp(2m,q)$.
	Then $G_0$ acts trivially on $\bfZ(q^{1+2m}_+)$ by Lemma~\ref{lem:autnmq}.
	Thus, $G_0$ is a subgroup of $\Aut(N)^V$ as $U^{G_0}=U$.
	Hence, $\Aut(N)^V$ is transitive on non-zero vectors of $N/\bfZ(N)$.

	Let $H=\Aut(q^{1+2m}_+)^{\bfZ(q^{1+2m}_+)}\cong G/G_0\cong\GammaL(1,q)$.
	Then $H_U\leqs \Aut(N)^{\bfZ(N)}$.
	Since $H_U$ is transitive on non-zero vectors of $\bfZ(N)$, we have $\Aut(N)^{\bfZ(N)}$ is transitive on non-identity elements of $\bfZ(N)$.
	By Lemma~\ref{lem:spext}\,(iv), we obtain that $N$ is a $3$-orbit group.

	Since $N$ is a $3$-orbit group, $\Aut(N)^{V}$ has a normal subgroup $G_1\cong\Sp(2m',q_0)$ with $q_0=p^{n'}$ such that $m'n'=mn$ by Lemma~\ref{lem:spext}\,(iv) and Theorem~\ref{thm:transmod}.
	Since $G_0\leqs \Aut(N)^{V}$, we have $n'\leqs n$ and $G_0\leqs G_1$.
	Suppose that $n'> n$.
	By Lemma~\ref{lem:nmpncgamma}\,(i), $N\cong (q_0)^{1+2m'}_+/U_1$ for some $U_1<\bfZ((q_0)^{1+2m'}_+)$ where $q_0=p^{n'}$.
	Then $(q_0)^{1_2m'}_+\cong q^{1+2m}_+/U_2$ with $U_2$ a subfield hyperplane of $\bfZ(q^{1+2m}_+)$ with respect to $\bbF_{q_0}$ by Theorem~\ref{thm:isos}\,(i).
	Hence, we have
	\[q^{1+2m}_+/U\cong N\cong (q_0)^{1+2m'}_+/U_1\cong (q^{1+2m}_+/U_2)/U_1.\]
	We may assume that the preimage of $U_1$ in $\bfZ(q^{1+2m}_+)$ is $U$.
	Then $U_2<U$ is a subfield hyperplane, a contradiction.
	Thus, we have $n=n'$ and $G_0=G_1\cong\Sp(2m,q)$.
	Therefore, $\Aut(N)^V\cong G_0{:}\GammaL(1,q)_U$, and then the lemma holds by Lemma~\ref{lem:spext}\,(i).
\end{proof}

Now we are in a position to prove Theorem~\ref{thm:3orbits}.
\begin{proof}[Proof of Theorem~\ref{thm:3orbits}]
	Let $N$ be a finite $3$-orbit group.
	Note that the set of orders of elements of $N$ can only be $\{1,p,q\}$, $\{1,p,p^2\}$ and $\{1,p\}$ for some primes $p$ and $q$.
	Hence, $N$ is either a $\{p,q\}$-group or a $p$-group.
	
	If $N$ is a $\{p,q\}$-group, then Lemma~\ref{lem:pq} shows that $N$ is in the first row of Table~\ref{tab:3aogp}.
	If $N$ is an abelian $p$-group, then $N$ is in the second row by Lemma~\ref{lem:ab}.
	If $N$ is a non-abelian $2$-group, then $N$ is one of the groups in rows~3-5 by Theorem~\ref{thm:2group}.

	Assume that $N$ is a non-abelian $p$-group for odd prime $p$.
	Then $N$ is a special $p$-group of exponent $p$ by Lemma~\ref{lem:3aop}.
	If $\Aut(N)$ is solvable, then $N\cong q^{1+2}_{+}$ by Theorem~\ref{thm:revdorn} and is in the 7-th row of Table~\ref{tab:3aogp}.
	If $\Aut(N)$ is non-solvable, then $G=\Aut(N)^{N/\bfZ(N)}$ is non-solvable by Lemma~\ref{lem:spext}\,(i).
	Lemma~\ref{lem:spext}\,(iii) shows that $G$ acts transitively on non-identity elements of both $N/\bfZ(N)$ and $\bfZ(N)$, and Lemma~\ref{lem:spext}\,(iv) shows that $\bfZ(N)$ is isomorphic to a quotient $\bbF_pG$-module of $\Lambda_{\bbF_p}^2(N/\bfZ(N))$.
	By Theorem~\ref{thm:transmod}, one of the following cases holds.
	\begin{enumerate}
		\item[(a)] $|\bfZ(N)|=1$;
		\item[(b)] $|\bfZ(N)|=|N/\bfZ(N)|$;
		\item[(c)] $G$ has a normal subgroup $G_0\cong \Sp(2m,p^n)$, and $|\bfZ(N)|\leqs p^n$.
	\end{enumerate}
	If~(a) holds, then $N$ is an extraspecial $p$-group of exponent $p$, which is in the 7-th row.
	If~(b) holds, then $N$ is in the 6-th row by Lemmas~\ref{lem:slunique} and~\ref{lem:autoapsl}.
	Assume that (c) holds.
	Then $N$ is a quotient of $q^{1+2m}_{+}/U$ for some $U\leqs\bfZ(q^{1+2m}_{+})$ by Lemma~\ref{lem:nmpncgamma}\,(ii), and $q^{1+2m}_{+}$ is listed in the 7-th row.
	If $U\neq 1$ and contains some subfield hyperplane $U_0$, then $N$ is a quotient of $q^{1+2m}_{+}/U_0\cong (q_0)^{(1+2mn_0)}_+$ by Theorem~\ref{thm:isos}\,(i).
	Thus, we may assume that $U$ contains no subfield hyperplanes.
	Therefore, Lemma~\ref{lem:quotaut} implies the last row of Table~\ref{tab:3aogp}, and we complete the proof.
\end{proof}

We now prove Corollary~\ref{coro:cgammaspfinal}.
\begin{proof}[Proof of Corollary~\ref{coro:cgammaspfinal}]
	Let $N_1=N/U$, and let $V=N/\bfZ(N)=N_1/\bfZ(N_1)$.
	Then $\Aut(N)^V\cong\CGammaSp(2m,q)$ by Lemma~\ref{lem:nmpncgamma}\,(iii).
	Since $\Sp(2m,q)\leqs \Aut(N)^V$ acts trivially on $\bfZ(N)$, it follows that $\Sp(2m,q)\leqs\Aut(N_1)^V$ and $\Aut(N_1)^V$ acts transitively on non-zero vectors of $V$.
	Thus, $N_1$ is $3$-orbit if and only if $\Aut(N_1)^{\bfZ(N_1)}$ is transitive on non-identity elements of $\bfZ(N_1)=\bfZ(N)/U$ by Lemma~\ref{lem:spext}\,(iii).
	It suffices to prove that $(H_U)^{\bfZ(N_1)}=\Aut(N_1)^{\bfZ(N_1)}$.

	Let $U_0$ be a subfield hyperplane with respect to subfield $\bbF_{p^d}$ for some $d\mid n$ such that $U_0\leqs U$ and $U_1=U/U_0$ has no subfield hyperplanes of $\bfZ(N)/U_0\cong \bbF_{p^d}^+$.
	Let $N_0=N/U_0$, and let $H_0=\Aut(N_0)^{\bfZ(N)/U_0}$.
	By Theorem~\ref{thm:isos}\,(i) and Lemma~\ref{lem:nmpncgamma}\,(iii), we have 
	\[N_0\cong (q_0)^{1+2mn/d}_+\mbox{ and }H_0\cong\GammaL(1,q_0).\]
	Note that $H_{U_0}^{\bfZ(N)/U_0}=H_0\cong\GammaL(1,q_0)$ by Lemma~\ref{lem:subequal}.
	Hence, we have 
	\[(H_U)^{\bfZ(N_1)}=(H_U)^{\bfZ(N)/U}=(H_0)_{U_1}^{\bfZ(N_0)/U_1}=(H_0)_{U_1}^{\bfZ(N_1)}.\]
	Since $U_1$ has no subfield hyperplanes, Lemma~\ref{lem:quotaut} implies that 
	\[\Aut(N_1)^{\bfZ(N_1)}= ((H_0)_{U_1})^{\bfZ(N_1)}=(H_U)^{\bfZ(N_1)}.\]
	Therefore, $N/U$ is $3$-orbit if and only if $H_U$ is transitive on $\bfZ(N)/U$.
\end{proof}

At last of this paper, we construct two families of groups in the 8-th row of Table~\ref{tab:3aogp}.

Let $V=\bbF_{p^n}^+$, and let $1\neq U<V$.
Note that $U$ is a subfield hyperplane when $n=2$.
For $n=3$, $U$ contains no subfield hyperplanes if and only if $\dim U=1$.
If $\dim U=1$, then $\GL(1,p^3)_U\cong\bbZ_{p-1}$, and hence $|\GammaL(1,p^3)_U|\leqs 3(p-1)$.
Notice that $|V/U|-1=p^2-1>3(p-1)$ in this case.
Therefore, smallest groups $q^{1+2m}_+/U$ in Example~\ref{exam:quot} satisfies $q\geqs 3^4$. 

When $V=\bbF_{3^4}^+$, let $\lambda$ be a generator of $\bbF_{3^4}^\times$ with minimal polynomial $\lambda^4-\lambda^3-1=0$.
Assume that $x$ is the operator of multiplying $\lambda$.
Then $V$ is a direct sum of the following two subspaces:
\[U_1=\langle 1+\lambda^2,2\lambda+\lambda^2+2\lambda^3\rangle\mbox{ and }U_2=\langle 2+\lambda,1+\lambda+\lambda^2\rangle.\]
We remark that each of $U_1$ and $U_2$ contains no subfield hyperplanes, and $H_{U_i}\cong\bbZ_8$ acts transitively on non-zero vectors of $V/U_i$ for $i=1,2$ where $H=\GammaL(1,3^4)$.
Hence, $U_1$ and $U_2$ admit two smallest $3$-orbit groups (of order $3^{10}$) in the 8-th row of Table~\ref{tab:3aogp}, which are $3^{1+2}_+/U_1$ and $3^{1+2}_+/U_2$, respectively.
Notice that there exists $h\in N_{H}(L)$ such that $U_1^h=U_2$, the above two groups are isomorphic.
We present a generating relation for the above group below:
\[\begin{aligned}
	N=\langle &x_i,y_j,z_k\mbox{ for $1\leqs i,j,k\leqs 4$ and $1\leqs k\leqs 4$}\mid x_i^3=y_i^3=z_i^3=1\mbox{ and }\\
	&[x_i,x_j]=[y_i,y_j]=[z_i,x_j]=[z_i,y_j]=1\mbox{ for each $1\leqs i,j\leqs 4$;}\\
	&[x_1,y_1]=z_1,\ [x_1,y_2]=[x_2,y_1]=z_2,\ [x_1,y_3]=[x_2,y_2]=[x_3,y_1]=z_3,\\
	&[x_1,y_4]=[x_2,y_3]=[x_3,y_2]=[x_4,y_1]=z_4,\ [x_2,y_4]=[x_3,y_3]=[x_4,y_2]=z_1z_4,\\
	&[x_3,y_4]=[x_4,y_3]=z_1z_2z_4,\ [x_4,y_4]=z_1z_2z_3z_4,\ z_1z_3=z_2^2z_3z_4^2=1\rangle.
\end{aligned}\]
We thank to Professor O'Brien, who told us that the \textit{TameGenus package} in \Magma\  provides some functions dealing with $p$-groups having nilpotency class $2$, exponent $p$ and center of order $p^2$.
The readers can use the \textbf{TGAutomorphismGroup} command in \Magma\ to check that $\Aut(N)\cong 3^{2\times 8}{:}(\Sp(2,3^4){:}\bbZ_8)$.

We now construct another family of groups in the 8-th row of Table~\ref{tab:3aogp} by obtaining a family of desired subspaces $U$ of $\bbF_q^+$.

\begin{example}\label{exam:nonsubfield}
	Let $p$ and $r$ be distinct odd primes such that $r\nmid p-1$, and let $q=p^n$ with $n=p^r-1$.
	Set $V=\bbF_{q}^+$, and let $H=\GammaL(1,q)$ acting naturally on $V$.
	Let $y\in H$ admitted by the Frobenius automorphism in $\Aut(\bbF_q)$.
	Then there exists an $r$-dimensional faithful $\bbF_p\langle y\rangle$-submodule $R$ of $V$.
	Let $U$ be an $\bbF_p L$-complement of $R$ in $V$.
\end{example}

With definitions in the above example, we observe that $|V/U|-1=|R|-1=p^r-1=|y|$.
Hence, $\langle y\rangle$ acts transitively on non-zero vectors of $V/U$, and then $q^{1+2m}_+/U$ is a $3$-orbit group for any positive integer $m$ by Corollary~\ref{coro:cgammaspfinal}.
Since $\dim R=r$ is not a divisor of $\dim V=p^r-1$, $q^{1+2m}_+/U$ lies in 8-th row of Table~\ref{tab:3aogp}.
Remark that the minimum integer $n$ satisfying the conditions in Example~\ref{exam:nonsubfield} is $n=p^r-1=5^3-1=124$.
Therefore, the smallest $3$-orbit group admitted by Example~\ref{exam:nonsubfield} has order $5^{251}$ with center of order $5^3$.


\end{document}